\documentclass[a4paper,oneside,10pt]{article}
\usepackage[english]{babel}
\usepackage[latin1]{inputenc}
\usepackage[T1]{fontenc}
\usepackage{amsmath}
\usepackage{amsfonts}
\usepackage{dsfont}
\usepackage{amssymb}
\title{Dirichlet-like space and capacity in complex analysis in several variables}
\usepackage{dsfont}
\author{Gabriel Vigny}
\begin{document}
\newtheorem{theorem}{Theorem}[section]
\newtheorem{proposition}[theorem]{Proposition}
\newtheorem{defi}[theorem]{Definition}
\newtheorem{corollaire}[theorem]{Corollary}
\newtheorem{lemme}[theorem]{Lemma}
\newtheorem{Remark}[theorem]{Remark}
\maketitle

\begin{abstract}
 For a Kähler manifold $X$, we study a space of test functions $W^*$ which is a complex version of $W^{1,2}$. We prove for $W^*$ the classical results of the theory of Dirichlet spaces: the functions in $W^*$ are defined up to a pluripolar set and the functional capacity associated to $W^*$ tests the pluripolar sets. This functional capacity is a Choquet capacity. \\
 The space $W^*$ is not reflexive and the smooth functions are not dense in it for the strong topology. So the classical tools of potential theory do not apply here. We use instead pluripotential theory and Dirichlet spaces associated to a current.
\end{abstract}


\section{Introduction}
The theory of Dirichlet spaces has been developped in the 1960's as a powerful tool in potential theory (see e.g \cite{Deny}). Its developments led to the theory of functional capacity for the Sobolev spaces and to the theory of pointwise value for functions in $W^{k,p}$ (see \cite{FZ1}, \cite{Zie} and \cite{AH1}). Those topics are useful in functional analysis, PDE and dynamics. In complex analysis in several variables, the work of Bedford and Taylor (\cite{BT1}) is a non-linear generalization of the Newtonian capacity in potential theory where the Laplacian is replaced by the Monge-Ampère operator $(dd^c)^k$ and the subharmonic functions by the plurisubharmonic functions (psh for short). Since then, the Bedford-Taylor capacity has been greatly used and studied. Capacities have become a classical tool in complex analysis and dynamics in several variables. Several capacities have been introduced since then (see \cite{SW1}, \cite {Ale}, \cite{DS2} and  \cite{GZ1}). Still, none of those capacities are "functional capacities" as in \cite{FZ1}. The main purpose of this article is to introduce such a capacity for compact Kähler manifolds.\\

For that, we present here a complex Sobolev space $W^*$ introduced by Dinh-Sibony in \cite{DS4}. Their purpose was to give a new proof of the decay of correlations for the dynamics of meromorphic maps. Let $X$ be a Kähler manifold and $\omega$ a Kähler form on $X$. The space $W^*$ is the subspace of $W^{1,2}$ of functions $\varphi$ such that there exists a positive closed current of finite mass $T_\varphi$ satisfying:
$$i\partial \varphi \wedge \bar{\partial} \varphi \leq T_\varphi,$$
and we define on $W^*$ the norm:
$$ \|\varphi \|_{*}^2=  \|\varphi \|_{L^2}^2 + \inf \Big\{ m(T)| T \ \text{closed, satisfying} \ i\partial \varphi \wedge \bar{\partial} \varphi \leq T \Big\},$$ 
where $m(T):=\int_X T\wedge \omega^{k-1}$. This definition seems more fitted to the complex case because it is independent of holomorphic change of coordinates. Our purpose in this article is to adapt the theory of Dirichlet spaces to the space $W^*$. In particular we want to show that the space $W^*$ characterizes pluripolar sets just like the space $W^{1,2}$ characterizes polar sets. Namely, we show that functions in $W^*$, a priori defined almost everywhere, can in fact be defined up to a pluripolar set, and that the functional capacity $C$ associated to $W^*$ defined for a Borel set $E$ by:
$$C(E)=\inf \Big\{ \|\varphi \|^2_*, \  \varphi \leq -1 \ \text{ a.e on some neighborhood of} \ E, \ \varphi \leq 0 \ \text{on} \ X\Big\}$$
satisfies $C(E)=0$ if and only if $E$ is pluripolar. On the other hand, they are some main differences with the classical case: smooth functions are not dense for the strong topology in $W^*$, and we will see in corollary \ref{nonreflexive}, that the space $W^*$ is not reflexive. So all the classical proofs in the theory of functional capacities fail and we will need to use other tools, especially pluripotential theory, Bedford-Taylor capacity and Dirichlet spaces associated to a current (see \cite{Oka1}). \\

The space $W^*$ appears as a space of test functions suited to complex analysis and dynamics: it is in a way very similar to $W^{1,2}(\Sigma)$ where $\Sigma$ is a Riemann surface, so we will stress on the similarities. Let us now detail the strucure of the article.

In  section \ref{section 2}, we study the basic properties of the space $W^*$ in both the local and compact case. In particular, we show that it is a Banach space and we introduce a notion of weak convergence for which compactness results hold. Then we show that the elements of $W^*$ are in $BMO$, so they are in $L^p$ for all $p<+\infty$. We give examples showing that the elements of $W^*$ are not in $VMO$ in general, that smooth functions are not dense in $W^*$ for the strong topology and that $W^*$ is not reflexive. We conclude by a theorem of weak density of smooth functions in $W^*$ for a compact Kähler manifold.

We consider the local situation in section \ref{quasi-continuity}. We prove two of our main results: the functions in $W^*$ are in the Dirichlet spaces associated to some positive closed currents and are quasi-continuous for the Bedford-Taylor capacity. This allows us to define for each function in $W^*$ a value at every point outside a pluripolar set. In fact, functions in $W^*$ are pluri-finely continuous outside a pluripolar set, which leads to interesting properties for the size of their Lebesgue set.

In section 4, we consider the case where $X$ is a compact Kähler manifold. We develop the study of the functional capacity $C$ for $W^*$ for which pluripolar sets are the sets of zero capacity. We show that it defines a Choquet capacity. The continuity result $C(\cup E_i)=\lim C(E_i)$ for $(E_i)$ an increasing sequence of Borel sets is the main difficulty, and it uses the results of the previous section. We show that this capacity and the Bedford Taylor capacity are comparable. We briefly explain how to extend these results to the local case. Finally, we introduce its dual capacity using classical arguments. 

In an appendix, we recall essential facts of the theory of Dirichlet spaces in the setting of Dirichlet space associated to a positive closed current (see \cite{Deny} and particularly \cite{Oka1}). We give all the proofs for the reader's convenience and we stress on what would make them fail in the case of $W^*$. More precisely, the classical approach consists in constructing the capacity first and using it to prove quasi-continuity results whereas our approach for $W^*$ is to prove quasi-continuity results first. The reader not familiar with Dirichlet spaces may start with this appendix.\\

\noindent {\bf Acknowledgements.} The author thanks Tien-Cuong Dinh and Nessim Sibony for their advices during the preparation of this article.

\section{General setting}	\label{section 2}
\subsection{Definitions and first results}
	Let $X$ be a Kähler manifold of dimension $k$ which is either a bounded open set of $\mathbb{C}^k$ or a compact manifold. We denote by $d^c$ the operator $d^c:=i/(2\pi)(\bar{\partial}-\partial)$.	We let $\omega$ be a Kähler form on $X$ such that $\int_X \omega^k=1$ (for a bounded domain in $\mathbb{C}^k$, we use the Kähler form $\omega=c dd^c \|z\|^2$, $c>0$). Let $W^{1,2}$ be the Sobolev space of real functions in $L^2$ whose differential in the sense of currents is a form with $L^2$ coefficients. We define the norm $\| \varphi \|^2_{W^{1,2}}=\| \varphi \|_{L^2}^2+ \|d\varphi \|_{L^2}^2$  on $W^{1,2}$. Define $W^{*}$ as the subspace of $W^{1,2}$ corresponding to the functions $\varphi \in W^{1,2}$ such that there is a bidegree (1,1) closed current $T$ of finite mass on $X$ with:
	\begin{equation}
	 \label{definition}
	 i\partial \varphi \wedge \bar{\partial} \varphi \leq T
	\end{equation}
in the sense of currents \cite{DS4}. Observe that $T$ satisfying (\ref{definition}) is positive since the left-hand side is positive and if $X$ is compact it is always of finite mass. Observe that when $k=1$, $W^*=W^{1,2}$ since every $(1,1)$ form is closed in dimension 1. The set of currents satisfying (\ref{definition}) is in fact convex and closed in the sense of currents since a weak limit of a sequence of positive (resp. closed) currents is positive (resp. closed). 

For $\varphi \in W^{*}$, we use the notation: 
	$$ \|\varphi \|^2_{*}=  \|\varphi \|_{L^2}^2 + \inf \Big\{ m(T)| T \ \text{closed, satisfying (\ref{definition})} \Big\}$$ 
	where $m(T):=\int_X T\wedge \omega^{k-1}$ is the mass of $T$. Observe that the infimum in the definition of $\|\varphi\|_{*}$ is reached because it is taken over a closed set. Furthemore, a current minimal in (\ref{definition}) will not charge hypersurfaces. Indeed by Siu's theorem \cite{Siu}, it can then be written $T_1+T_2$ where the $T_i$ are positive closed currents with $T_1$ a current of integration on a countable union of hypersurfaces and $T_2$ giving no mass to hypersurfaces, and $T_2$ will satisfies (\ref{definition}). Clearly, there exists a constant $A>0$ such that for $\varphi \in W^{*}$, we have $ \|\varphi \|_{W^{1,2}} \leq  A \|\varphi \|_{*}$. We have the proposition:
\begin{proposition}\label{norm}
The function $\varphi \mapsto  \|\varphi \|_{*}$ is a norm on $W^{*}$ and $W^{*}$ is complete with respect to this norm.

 \end{proposition}
	\emph{Proof.} For the first assertion, only the triangular inequality has to be proved. Let $\varphi$ and $\psi$ in $W^{*}$, and  $T_{\varphi}$ and $T_{\psi}$ be minimal for the mass in (\ref{definition}), then:
\begin{eqnarray*}
	  i\partial (\varphi +\psi) \wedge \bar{\partial}( \varphi +\psi) = i\partial \varphi \wedge \bar{\partial} \varphi + i\partial \psi \wedge \bar{\partial} \psi + i(\partial \varphi \wedge \bar{\partial} \psi + i\partial \psi \wedge \bar{\partial} \varphi)
	  \end{eqnarray*}
	If $T_{\varphi}$ is zero, $\varphi$ is constant and the result is clear. Otherwise take $c=  \sqrt{\frac{m(T_{\psi})}{m(T_{\varphi})}}$. By Cauchy-Schwarz inequality:
	$$i(\partial \varphi \wedge \bar{\partial} \psi + i\partial \psi \wedge \bar{\partial} \varphi) \leq c i\partial \varphi \wedge \bar{\partial} \varphi + \frac{1}{c} i\partial \psi \wedge \bar{\partial} \psi.$$ 
Hence:
	 \begin{eqnarray*}
	i\partial (\varphi +\psi) \wedge \bar{\partial}( \varphi +\psi) \leq (1+c) T_{\varphi} +(1+ \frac{1}{c}) T_{\psi}. 
	 \end{eqnarray*}
The left-hand side is a positive closed current of mass $(\sqrt{m(T_{\varphi})}+\sqrt{m(T_{\psi})})^2$ which gives the result. \\
	
	For the second assertion, observe that if $(\psi_n)_{n\in \mathbb{N}}$ is a Cauchy sequence in $W^{*}$, it is a Cauchy sequence in $W^{1,2}$ which is complete, so it converges in $W^{1,2}$ to a function $\psi \in W^{1,2}$. For every $\varepsilon>0$, there is an integer $N$ such that for $n$ and $m$ larger than $N$ we have $i\partial (\psi_n-\psi_m) \wedge \bar{\partial}(\psi_n-\psi_m) \leq T_{n,m}$ where $T_{n,m}$ is a closed current of mass less than $\varepsilon$. We let $n$ go to infinity. Since $(\psi_n-\psi_m)_n$ converges in $W^{1,2}$ to $\psi-\psi_m$, we have that $(i\partial (\psi_n-\psi_m) \wedge \bar{\partial}(\psi_n-\psi_m))_n$ converges in $L^1$ thus in the sense of currents to $i\partial (\psi-\psi_m) \wedge \bar{\partial}(\psi-\psi_m)$. And we can extract a subsequence of  $(T_{n,m})_n$ which converge in the sense of currents to some closed current $T_m$ of mass less than $\varepsilon$ satisfying (\ref{definition}) for $\psi-\psi_m$ since a weak limit of positive current is positive. This gives that $\psi$ is in fact in $W^*$ and that $(\psi_m)$ converges to $\psi$ in $W^*$.
	$\Box$ \\
	
The following result is deduced from the previous proof:
	\begin{corollaire}\label{monotone}
If $(\varphi_n)$ is a bounded sequence in $W^*$ converging in $W^{1,2}$, then its limit is in $W^*$.
  \end{corollaire}
We will see in section \ref{examples} that smooth functions are not dense in $W^*$ and natural sequences do not converge for this topology. So we will use the following weak convergence:
  \begin{defi}
 Let $(f_n)$ be a sequence in $W^*$ and $f \in W^*$, we write $f_n \rightharpoondown f$ if  $(f_n)$ converges weakly to $f$ in $W^{1,2}$ and $( \|f_n \|_*)$ is a bounded sequence.
 \end{defi}
Of course, it is the same to ask for $(f_n)$ to converge in the sense of distributions and for $( \|f_n \|_*)$ to be a bounded sequence, but our definition is more convenient. The previous definition is interesting because of the following compactness result:
 \begin{proposition}
 Let $(f_n)$ be a bounded sequence in $W^{*}$. Then there exist a subsequence $(f_{n_j})$ and $f\in W^{*}$ such that $f_{n_j}\rightharpoondown f$. Furthermore, we have $\|f \|_* \leq \liminf \|f_{n_j} \|_*$.
 \end{proposition}
\emph{Proof.} Since $(f_n)$ is bounded in $W^{*}$, it is also bounded in $W^{1,2}$. Taking a subsequence, we can assume that $(f_n)$ converges weakly in $W^{1,2}$ to $f \in W^{1,2}$. Let $T_n$ be a closed current satisfying (\ref{definition}) of minimal mass for $f_n$. We can assume that $(T_n)$ and  $(i\partial  f_n\wedge \bar{\partial} f_n)$ converge in the sense of currents to some limits $T$ and $\Theta$ with $T$ positive and closed and $\Theta$ positive. Let $\Psi$ be a positive test form of bidegree $(k-1,k-1)$, we want to show that $\langle i\partial  f\wedge \bar{\partial} f,\Psi \rangle \leq \langle \Theta,\Psi \rangle$, which will conclude the proof since $\Theta \leq T$.

 By the definition of positive forms: $ \langle i\partial  (f-f_n)\wedge \bar{\partial} (f-f_n),\Psi \rangle \geq 0$, we expand:
 $$ \langle i\partial f\wedge \bar{\partial} f_n+i\partial f_n\wedge \bar{\partial}f ,\Psi \rangle \leq  \langle i\partial f\wedge \bar{\partial} f,\Psi \rangle+\langle i\partial f_n\wedge \bar{\partial} f_n,\Psi \rangle.$$
We have that $(\langle i\partial f\wedge \bar{\partial} f_n,\Psi \rangle)$ goes to $\langle i\partial f\wedge \bar{\partial} f,\Psi \rangle$ because $\partial f\wedge \Psi$ has coefficients in $L^2$ and $(\bar{\partial} f_n)$ has coefficients weakly converging in $L^2$. Similarly $(\langle i\partial f_n\wedge \bar{\partial} f,\Psi \rangle)$ goes to $\langle i\partial f\wedge \bar{\partial} f,\Psi \rangle$. Letting $n\to \infty$ gives:
 $$ \langle i\partial f\wedge \bar{\partial} f+i\partial f\wedge \bar{\partial}f ,\Psi \rangle \leq \langle i\partial f\wedge \bar{\partial} f,\Psi \rangle+\langle\Theta,\Psi \rangle $$ 
 which concludes the proof. $\Box$\\
 
Let $U$ be an open set in $\mathbb{C}^k$, $U'\Subset U$, and $\varphi \in W^*(U)$. Take $\chi$ a non negative smooth radial function with compact support in $\mathbb{C}^k$ such that $\int \chi =1$. Define $\chi_\varepsilon(z):=\varepsilon^{-2k}\chi(z/ \varepsilon)$ for $\varepsilon>0$ and put $\varphi_\varepsilon=\varphi * \chi_\varepsilon$ (well defined in $U'$ for $\varepsilon$ small enough), then $\varphi_\varepsilon$ is smooth. Let $(\varepsilon_n)$ be a sequence decreasing to zero and define $\varphi_n=\varphi_{\varepsilon_n}$. It is classical that $(\varphi_n)$ converges to $\varphi$ in $W^{1,2}(U')$. Let $T$ be such that $i\partial \varphi \wedge \bar{\partial}\varphi \leq T$ and let $v$ be a psh potential of $T$ (that is $i\partial \bar{\partial} v =T$), we define $T_n=T *\chi_{\varepsilon_n}$ and $v_n =v *\chi_{\varepsilon_n}$ so that $i\partial \bar{\partial} v_n =T_n$. Then $(v_n)$ decreases to $v$ and $(T_n)$ converges to $T$ in the sense of currents. In particular, 
$$\int_{U'} T\wedge \omega^{k-1}\leq \lim \int_{U'} T_n\wedge \omega^{k-1}\leq\int_{U} T\wedge \omega^{k-1}$$
 Using the previous notations, we can now state a regularization lemma: 
\begin{lemme}\label{regularization}
\begin{enumerate}
\item Let $U$ be an open set in $\mathbb{C}^k$. Then for $U'\Subset U$, and $\varphi \in W^*(U)$, the sequence of smooth functions $(\varphi_n)$ converges weakly to $\varphi$ in $W^*(U')$. Furthermore, we have that $i\partial \varphi_n \wedge \bar{\partial}\varphi_n \leq T_n$. In particular, we have the inequalities $\|\varphi \|_{W^*(U')} \leq \lim \| \varphi_n\|_{W^*(U')} \leq \|\varphi \|_{W^*(U)} $. Finally, we can choose the potential $v_n$ of the currents $T_n$ so that $(v_n)$ decreases to the potential $v$ of $T$. 
\item Consider the projective space $\mathbb{P}^k$. Let $\varphi \in W^*(\mathbb{P}^k)$, then there exists a sequence of smooth functions $(\varphi_n)$ converging weakly to $\varphi$ in $W^*(\mathbb{P}^k)$ and $\lim \| \varphi_n \|_* = \| \varphi \|_*$.
\end{enumerate}
\end{lemme}
\emph{Proof.} For the first case, we have seen in the proof of proposition \ref{norm} that if $i\partial \varphi \wedge \bar{\partial}\varphi \leq T_\varphi$ and $i\partial \psi \wedge \bar{\partial}\psi \leq T_\psi$ then:
$$i\partial (\frac{\varphi+\psi}{2}) \wedge \bar{\partial}(\frac{\varphi+\psi}{2}) \leq (\frac{T_\varphi+T_\psi}{2}).$$
Approximating $\chi_{\varepsilon_n}$ by a finite sum and using that convexity property, we get $i\partial \varphi_n \wedge \bar{\partial}\varphi_n \leq T_n$. The rest follows.

In the second case, we apply the same argument with an approximation of $\varphi$, using an approximation of the identity in $\text{Aut}(\mathbb{P}^k)$. The current $T_n$ satisfying $i\partial \varphi_n \wedge \bar{\partial}\varphi_n \leq T_n$ converges to $T$ (as an average of the composition of $T$ by automorphisms of $\mathbb{P}^k$) hence $m(T_n)=\langle T_n, \omega^{k-1} \rangle \to \langle T_n, \omega^{k-1} \rangle = m(T)$. So $\lim \|\varphi_n\|_*=\|\varphi\|_*$. $\Box$\\

In the first case, we cannot expect in general the equality $\lim \| \varphi_n \|_{W^*(U')} = \| \varphi \|_{W^*(U)}$ since the best current in $U'$ for $\varphi$ is not necessarily the restriction of the best current in $U$ (take a non constant function on $U$ that vanishes on $U'$). We will prove a version of the above result in the case of compact Kähler manifold in section \ref{weak_density}.\\ 

\subsection{Functions in  $W^*$ and $BMO$} 
For $x \in \mathbb{C}^k$ and $T$ a positive closed (1,1)-current defined on some neighborhood of $U$, let $\nu(T,x,r):=r^{2(1-k)}\int_{B(x,r)} T\wedge (dd^c\|z\|^2)^{k-1}$ where $B(x,r)$ is the ball of center $x$ and of radius $r$. We know the quantity $\nu(T,x,r)$ decreases to the \emph{Lelong number} of $T$ at $x$ when $r$ decreases to $0$ (see \cite{dem2} for details).\\
 
 Let $U$ be some bounded open set in $\mathbb{R}^n$ with a riemannian metric $g$. Let $B$ be a ball contained in $U$, $|B|$ its volume and $f \in L^1(U)$ a real-valued function. We write $m_B(f)=\frac{1}{|B|}\int_B f$ the mean value of $f$ in the ball $B$. By definition, a function $f\in L^1(U)$ is in $BMO(U)$ (bounded mean oscillation) if there exists a constant $C_f$ such that for any ball $B(x,r)$ contained in $U$, we have that:
  $$\frac{1}{|B|}\int_B |f-m_B(f)| \leq C_f. $$
  We denote by $\|f\|_{BMO}$ the infimum of the numbers $C_f$. It is known that $\|f\|_{BMO}$ defines a semi-norm and that if  $f\in BMO(U)$ then there exists a constant $C'_f>0$ such that $e^{C'_f|f|}$ is in $L^1(U)$. More precisely, there exists a constant $b>0$ that depends only on $n$ such that for all $b'< b/\|f\|_{BMO}$, $e^{b'|f|}$ is in $L^1(U)$. In particular $f \in L^p(U)$ for all $p<\infty$ \cite{JN1}. Observe that $BMO(U)$ does not depend on the choice of $g$ and that we can extend the notion of $BMO$ to manifolds. We have the following proposition:
  \begin{proposition}
  Let $\varphi$ be in $W^*$. Then $\varphi$ is in $BMO$, consequently, $\varphi$ is in $L^p$ for all $p <\infty$   
  \end{proposition}  
  \emph{Proof.} Recall first that for a function in $W^{1,2}(U)$ where $U$ is an open set of $\mathbb{R}^n$, and for any ball $B\subset U$, we have the following Poincaré-Sobolev inequality (e.g. \cite{AH1}):
  \begin{eqnarray}\label{PS}
  \frac{1}{|B|}\int_B |\varphi-m_B(\varphi)| \leq C \frac{1}{|B|^{\frac{1}{2}-\frac{1}{n}}} \Big(\int_B \| \text{d}\varphi \|^2 \Big)^{\frac{1}{2}},
  \end{eqnarray}  
where $C$ is a constant that depends only on $n$ (in our case, $n=2k$). Using a covering if necessary, we can suppose that we are in an open set of $\mathbb{C}^k$. For $\varphi \in W^*$, $T_\varphi$ satisfying $i\partial \varphi \wedge \bar{\partial} \varphi \leq T_\varphi$, and $B$ a ball centered at $x$ of radius $r$, we deduce from (\ref{PS}) that:
  \begin{eqnarray}\label{bmo-lelong}
  \frac{1}{|B|}\int_B |\varphi-m_B(\varphi)| \leq C r^{1-k} \Big(\int_B T_\varphi \wedge (dd^c\|z\|^2)^{k-1}\Big)^{\frac{1}{2}}\leq C\sqrt{\nu(T_\varphi,x,r)}, 
 \end{eqnarray}
 and we know the quantity $\nu(T_\varphi,x,r)$ decreases to the Lelong number of $T_\varphi$ at $x$ when $r$ decreases to $0$. With the above notations, for any $C'_\varphi <b/\|\varphi\|_*$ where $b$ is a constant that depends only on $X$, then $e^{C'_\varphi|\varphi|}$ is in $L^1$.  $\Box$\\

 In particular, we see from the proof that if $T_\varphi$ has no positive Lelong number on $X$, then $\varphi$ is in fact $VMO$ (i.e: the mean oscillation is bounded and goes to zero when $r$ goes to zero). We will see that in general functions in $W^*$ are not in $VMO$ in the case of higher dimension. In dimension $1$, any function in $W^*=W^{1,2}$ is in $VMO$ (it is a consequence of the above proof).\\
  
\subsection{Examples, density and duality}\label{examples}
Lipschitz functions are in $W^*$. Furthermore, if $g:\mathbb{R}\to\mathbb{R}$ is Lipschitz and $f$ is in $W^*$ then $g\circ f\in W^*$. In particular we will use the fact that for $a\in \mathbb{R}$, $\max(f,a)$ is in $W^*$ with $\|\max(f,a) \|_* \leq \|f\|_*$ \cite{DS4}. 
For $f$ and $g$ smooth functions and $\varepsilon>0$, we let:
 \begin{eqnarray*}
 \text{max}_\varepsilon(f,g):= \frac{f+g+\sqrt{(f-g)^2+\varepsilon}}{2} \quad \text{and} \quad \text{min}_\varepsilon(f,g):= \frac{f+g-\sqrt{(f-g)^2+\varepsilon}}{2}.
\end{eqnarray*}
The functions $\text{max}_\varepsilon(f,g)$ and $\text{min}_\varepsilon(f,g)$ are smooth. A straightforward computation shows that:
$$i\partial (\text{max}_\varepsilon(f,g))\wedge \bar{\partial}(\text{max}_\varepsilon(f,g))+i\partial (\text{min}_\varepsilon(f,g))\wedge \bar{\partial}(\text{min}_\varepsilon(f,g))\leq i\partial f\wedge \bar{\partial}f+i\partial g\wedge \bar{\partial}g.$$
Letting $\varepsilon$ go to zero, we deduce that if $T_f$ and $T_g$ satisfy (\ref{definition}) for $f$ and $g$ then $T_f+T_g$ satisfies (\ref{definition}) for $\max(f,g)$ and $\min(f,g)$. The last property is local so by density it is true for any $f$ and $g$ in $W^*$. We deduce that $\max(f,g)$ and $\min(f,g)$ are in $W^*$ with the bound:
\begin{eqnarray}\label{extremum}
\|\max(f,g) \|_*^2\leq \|f\|_*^2+\|g\|_*^2 \quad \text{and} \quad \|\min(f,g) \|_*^2\leq \|f\|_*^2+\|g\|_*^2.
\end{eqnarray}

Now, if $\xi$ is a smooth function (even with compact support in the local case) and $f \in W^*$, then we do not know if $\xi f$ belongs to $W^*$ for $k\geq2$ in general. This is an important difference with the case of Sobolev spaces as partition of unity is a classical tool. Still, such rigidity is characteristic of complex analysis. We give now some less simple examples.  \\  

\noindent \emph{Example 1.} Let $X$ be either a compact Kähler manifold or a bounded open set in $\mathbb{C}^n$. Let $\varphi$ be a \emph{quasi plurisubharmonic} (qpsh for short) function on $X$, that is $\varphi$ is locally the sum of a psh function and a smooth function. So, we have $i\partial\bar{\partial}\varphi+C \omega \geq 0$ for some $C\geq 0$. If $\varphi$ is bounded (say $0\leq\varphi\leq 1$), then it is in $W^*$ because it satisfies $i\partial \varphi\wedge\bar{\partial}\varphi=i/2\partial\bar{\partial}(\varphi^2)-i\varphi\partial\bar{\partial}\varphi$ and the right-hand side term is bounded by the positive closed current $i/2\partial\bar{\partial}(\varphi^2)+C\omega$. 

Take now $\varphi \leq -1$ not necessarily bounded, then the  function $\psi= -\log (- \varphi)$ is in $W^{*}$. Indeed, since $i\partial \bar{\partial} \varphi \geq - C\omega$, we have:
\begin{eqnarray*}
 i \partial \psi \wedge \bar{\partial} \psi & = & \frac{i\partial \varphi \wedge \bar{\partial} \varphi }{|\varphi|^2}\\
 i\partial \bar{\partial} \psi              & = & -\frac{i\partial \bar{\partial} \varphi}{\varphi} + \frac{i\partial \varphi \wedge \bar{\partial} \varphi }{|\varphi|^2}.
\end{eqnarray*}
We have that $ i \partial \psi \wedge \bar{\partial} \psi =  i\partial \bar{\partial} \psi+ i\partial \bar{\partial} \varphi / \varphi\leq  i\partial \bar{\partial} \psi+C\omega$ which is of mass $C$. \\

\noindent \emph{Example 2.} Consider a bounded open set $X$  contained in the unit ball of $\mathbb{C}^n$ with the canonical Kähler form. Let $\varphi$ be the function defined by $(-\log |z_1|^2)^{\alpha}$ for $\alpha<1/2$. Then, $\varphi$ is in $W^*$ since it is in $L^2$ with:
$$i\partial \varphi \wedge \bar{\partial}\varphi= \frac{idz_1\wedge d\bar{z}_1}{|z_1|^2(-\log |z_1|^2)^{2-2\alpha}},$$
 which is closed and in $L^1$ if and only if $2-2\alpha>1$.\\
 
\noindent \emph{Example 3.} Consider a Kähler manifold $X$ of dimension 2 with some point $0\in X$. We consider the blow-up $\widetilde{X}$ of $X$ at $0$ and we denote by $\pi:\widetilde{X}\to X$ the standard projection. Let $H:=\pi^{-1}\{0\}$ denote the exceptional fiber. Let $\widetilde{f}$ be a smooth function on $\widetilde{X}$ (hence $\widetilde{f}\in W^*(\widetilde{X})$) so that $\int i\partial \widetilde{f} \wedge \bar{\partial}\widetilde{f}\wedge [H] >0$. Consider a current $\widetilde{T}$ of minimal mass satisfying (\ref{definition}) for $\widetilde{f}$.\\
Define $f:=\pi_* \widetilde{f}$ and $T:=\pi_* \widetilde{T}$, then $f\in W^*(X)$ since $i\partial f \wedge \bar{\partial}f \leq T$ indeed $i\partial f \wedge \bar{\partial}f$ gives no mass to 0. We want to compute the "Lelong number" $\lim_{r\to 0} \int_{\mathbb{B}_r} i\partial f \wedge \bar{\partial}f \wedge i \partial \bar{\partial} \log\|Z\|$, where $Z=(z,w)$ is a system of local coordinates. Recall that locally $\widetilde{X}$ is given as the set of points $\{((z,w),[x:y]), zy=wx\}$. In the chart where $x\neq 0$, we can write $\widetilde{X}$ as $\{((z,uz),[1:u])\}$. The projection $\pi$ in the $(z,u)$ coordinates is $(z,u)\mapsto (z,uz)$ , and $H$ is given by $z=0$. So we have:
\begin{eqnarray*}
 \int_{\mathbb{B}_r} i\partial f \wedge \bar{\partial}f \wedge i \partial \bar{\partial} \log\|Z\|&=&  \int_{\pi^{-1}(\mathbb{B}_r)} i\partial \widetilde{f} \wedge \bar{\partial}\widetilde{f} \wedge i \partial \bar{\partial} \log \| \pi(Z)\| \\
                                                                                                  &=&\int_{\pi^{-1}(\mathbb{B}_r)} i\partial \widetilde{f} \wedge \bar{\partial}\widetilde{f} \wedge i \partial \bar{\partial} \frac{1}{2}\log \big(|z|^2(1+|u|^2)\big)\\
                                                                                                  &\geq& \int_{\pi^{-1}(\mathbb{B}_r)} i\partial \widetilde{f} \wedge \bar{\partial}\widetilde{f} \wedge i \partial \bar{\partial} \log |z|.
 \end{eqnarray*}
 When $r$ goes to 0, the last integral goes to $\int i\partial \widetilde{f} \wedge \bar{\partial}\widetilde{f}\wedge [H]$ which is positive by our assumption. In particular any current satisfying (\ref{definition}) for $f$ has a Lelong number at zero. For example, take a smooth function $\widetilde{f}$ on $\widetilde{X}$ given by $\frac{|x|^2}{|x|^2+|y^2|}$ near $H$. It is smooth because it is given by $\frac{1}{1+|u|^2}$ in the chart where $x\neq 0$. Its push-forward $f$ is simply $\frac{|z|^2}{|z|^2+|w^2|}$, which is then in $W^*$. Recall the bound for the mean oscillation of $f$:
 $$\frac{1}{|B|}\int_B |f-m_B(f)| \leq C r^{1-k} \Big(\int_B T_f \wedge \omega^{k-1}\Big)^{\frac{1}{2}}.$$
  Since $f$ is homogeneous of order zero then $m_B(f)$ does not depend on $B$, so the function $|f-m_B(f)|$ is also homogeneous of order zero hence $m_B(|f-m_B(f)|)=A$ does not depend on $B$ and it is positive as $f-m_B(f)$ is not everywhere zero. So $f$ is not $VMO$. And if we apply the last formula to $f-f'$ where $f'$ is smooth, we see that the term $m_{B_r}(|f-f'-m_{B_r}(f-f'))|\geq m_{B_r}(|f-m_{B_r}(f)|-m_{B_r}(|f'-m_{B_r}(f')|))$ goes to $A$ when $r\to 0$. So a current $T_{f-f'}$ satisfying $i\partial (f-f')\wedge \bar{\partial}(f-f')\leq T_{f-f'}$ has a Lelong number bounded from below by a quantity which does not depend on $f'$, and hence a mass that does not depend on $f'$. The example is easily extended to higher dimension. So we proved the important proposition:
  \begin{proposition}
The space $W^*$ is not contained in $VMO$ and smooth functions are not dense in $W^*$ for the strong topology.
  \end{proposition}
The second assertion is also true for continuous functions in $W^*$ for the same reasons. We deduce the following fundamental corollary:
\begin{corollaire}\label{nonreflexive}
The space $W^*$ is not reflexive.
\end{corollaire}
\emph{Proof.} Assume it is reflexive. In this case, we consider the function $f$ above with support contained in some chart. So we are in the case of lemma \ref{regularization} and we take a sequence of regularizations $(f_n)$. This is a bounded sequence so we can extract a weakly (in the dual sense, not in our weak sense) converging sequence. Because it also converges in $W^{1,2}$, its limit is $f$. But since the weak closure and the strong closure of a convex set coincide, the limit $f$ would be in the strong closure of the smooth functions in $W^*$ which contradict the previous proposition. $\Box$

\begin{Remark} \rm Our weak topology is weaker than the dual weak topology, but it enjoys a compactness property so it is the right one to consider.  
\end{Remark}

\subsection{A density theorem}\label{weak_density}
\noindent We want to prove the (weak) density of smooth functions in $W^*$, the question was raised in \cite{DS4} where the authors advise to follow the arguments of \cite{DS5} which is what we do. So let us recall the results of \cite{DS5} we need first. For a compact Kähler manifold $X$ of dimension $k$, there exist two sequences of positive closed currents $(K_n^+)$ and $(K_n^-)$ of bidegree $(k,k)$ on $X \times X$ with coefficients in $L^1$ such that $(K_n^+-K_n^-)$ converges to the current of integration on the diagonal of $X\times X$. A precise description of the singularities of $K_n^{\pm}$ implies that for a positive closed current $T$ of any positive bidegree, the (positive closed) currents $T^{\pm}_n(x):= \int_{y\in X} K^{\pm}_n(x,y)\wedge T(y)$ satisfy:
$$T^+_n-T^-_n \to T.$$
Moreover, there exists a constant $c>0$ independent of $T$ and $n$ such that $m(T_n^{\pm})\leq c m(T)$.  

Furthermore, there exists $\delta>0$ such that if $T$ has measure coefficients then $T^{\pm}_n$ have coefficients in $L^{1+\delta}$, if $T$ has coefficients in $L^p$ then  $T^{\pm}_n$ have coefficients in $L^{q}$ (where $q=\infty $ if $p^{-1}+(1+\delta)^{-1}\leq 1$ and $p^{-1}+(1+\delta)^{-1}= 1+q^{-1}$ otherwise), if $T$ has coefficients in $L^\infty$ then $T^{\pm}_n$ are continuous forms, if $T$ is a continuous form then $T^{\pm}_n$ has $\mathcal{C}^1$ coefficients. Finally, the currents $K^{\pm}$ are smooth outside the diagonal of $X\times X$ and satisfy $\| K_n^{\pm}(.,y) \|_{L^1}\leq A$ where $A$ is a constant that does not depend on $n$ and $y$. 

In particular, for a function $f$ in $L^1$, let $(f^{\pm}_n)$ be the sequences defined by $f^{\pm}_n=\int_{y\in X} f(y)K^{\pm}_n(.,y)$. Then $(f_n):=(f^+_n-f^-_n)$ converges to $f$ in the sense of distributions. Indeed, if $f$ is continuous, the result is clear by weak convergence and one has the bound $\|f_n\|_{L^1} \leq A\|f\|_{L^1}$ so $f_n\to f$ in $L^1$.\\

\noindent Define $K_n:=K^{+}_n-K^{-}_n$ and let $\varphi\in W^*$ with $T_\varphi$ as in (\ref{definition}). Define 
$$\varphi_n:=\int_{y\in X} \varphi(y)K_n(.,y),$$ 
which is in $L^\infty$ since $\varphi$ is in $L^q$ for all $q\geq 1$, and $(\varphi_n)$ converges to $\varphi$ in the sense of distributions. Let $\pi_1$ and $\pi_2$ denote the canonical projections from $X\times X$ to each of its factor, then $\varphi_n=(\pi_1)^* (((\pi_2)_*\varphi) K_n)$. Since $K_n$ is closed and $\varphi\in W^{1,2}$, then $i\partial \varphi_n= (\pi_1)^* (((\pi_2)_* \partial \varphi )\wedge K_n)$ and $i\bar{\partial} \varphi_n= (\pi_1)^* (((\pi_2)_* \bar{\partial} \varphi )\wedge K_n)$. That is:
\begin{eqnarray*}
\partial \varphi_n= \int_{y\in X}  K_n(x,y)\wedge \partial \varphi(y) \quad \text{and} \quad \bar{\partial} \varphi_n= \int_{y\in X}  K_n(x,y)\wedge \bar{\partial} \varphi(y).
\end{eqnarray*}
We make the wedge product of this two terms, it is positive hence real so we can take the real part. We obtain a sum of terms of the form:
\begin{eqnarray*}
\int_{y,y'\in X} K^{\pm}_n(x,y)  \wedge K^{\pm}_n(x,y')\wedge \text{Re}\left( i\partial \varphi(y) \wedge \bar{\partial}\varphi(y')\right). 
                                      \end{eqnarray*}
Since the currents $K^{\pm}_n$ are positive, the last integral is less than:
\begin{eqnarray*}
\int_{y,y'\in X} K^{\pm}_n(x,y)  \wedge K^{\pm}_n(x,y')\wedge \frac{1}{2}\big(T_\varphi(y)+T_\varphi(y')\big) .
                                      \end{eqnarray*}
Since $\int_{y\in X} K_n(x,y)\leq A$, the integral is in turn less than:
$$ A \int_{y\in X} K_n^{\pm}(x,y)\wedge T_\varphi(y).$$
That integral defines a positive closed current, of mass controlled by the mass of $T_\varphi$ and with coefficients in $L^{1+\delta}$. In particular, we have that $i\partial \varphi_n \wedge \bar{\partial} \varphi_n$ is bounded by a positive closed current, of mass controlled by the mass of $T_\varphi$ and coefficients in $L^{1+\delta}$. We iterate the process, and we gain regularity until we get functions in $\mathcal{C}^1$ controlled by a current with coefficients in $\mathcal{C}^1$ (the number of iterations does not depend on $\varphi$). A small perturbation of the current $T_\varphi$ gives the following theorem:  
\begin{theorem}
Let $\varphi\in W^*(X)$. Then there exists a sequence of smooth functions $(\varphi_n)$ such that $\varphi_n \rightharpoondown \varphi$ in $W^*$. More precisely, there exists a constant $C_1$ that does not depend on $\varphi$ such that $\|\varphi_n\|_*\leq C_1 \|\varphi\|_*$. 
\end{theorem}

\section{Pointwise values}\label{quasi-continuity}
Let $U$ be an open pseudoconvex set of $\mathbb{C}^k$. The Bedford-Taylor capacity ${\rm cap}_{BT}$ (see \cite{BT1}) is defined by:
$${\rm cap}_{BT}(E):=\sup \Big\{ \int_E (dd^c u)^k | u \ \text{psh}, \ 0\leq u \leq 1 \Big\}$$
for $E \subset U$ a Borel set. It is subadditive. A set is of zero capacity if and only if it is pluripolar. Recall that a set $P$ is \emph{pluripolar} in $U$ if there is a psh function $v$ such that $P\subset\{v=-\infty\}$. And it is \emph{complete pluripolar} in $U$ if there is a psh function $v$ such that $P=\{v=-\infty\}$\\

A priori, a function $\varphi$ in $W^*$ is defined only almost everywhere. The purpose of this section is to show that if $\varphi$ is in $W^*$ then it is quasi-continuous on any open set $V\Subset U$ for the Bedford-Taylor capacity: there exists a representative $\widetilde{\varphi}$ of $\varphi$ such that $\forall \varepsilon >0$ there exists an open set $U_\varepsilon$ such that ${\rm cap}_{BT}(U_\varepsilon)\leq \varepsilon$ and $\widetilde{\varphi}$ restricted to $V \backslash U_\varepsilon$ is continuous. Moreover we show that two such representatives coincide outside a pluripolar set. 

We work locally, so let $U$ be a strongly pseudoconvex open set in $\mathbb{C}^k$. Let $0\leq \chi \leq 1$ be some smooth function with compact support in $U$. First we prove that $\chi\varphi$ can be seen as an element of the Dirichlet spaces associated to a class of positive closed currents. Then we prove the quasi-continuity. In particular, $\varphi$ is defined up to a pluripolar set. We will then prove a convergence lemma and a result on the Lebesgue set of $\varphi$. 

\subsection{Embedding in some Dirichlet spaces}\label{Embedding}

Let $S$ be a positive closed current of bidegree $(p,p)$ for $p\leq k-1$ in $U$ ($S$ is not necessarily of bounded mass). Denote by $H_S$ the completion of the smooth functions with compact support in $U$ for the hermitian norm $\|f\|_S^2:=\int_U i\partial f \wedge \bar{\partial}f \wedge S\wedge \omega^{k-p-1}$. It is a (real) Hilbert space. Let $u_1,\dots, u_{k-1}$ be bounded psh functions on $U$. Our purpose is to show that $\chi W^*$ can be embedded in the Dirichlet space $H_S$ when $S=i\partial \bar{\partial} u_1\wedge \ldots \wedge i\partial \bar{\partial} u_{k-1}$.  Let us stress that it is not clear that the quantity $i\partial \varphi \wedge \bar {\partial} \varphi \wedge S$ makes sense for $\varphi \in W^*$ since one cannot multiply a Radon measure by a form in $L^1$ in general. Let $U'$ be an open set of $U$ such that $\text{supp} (\chi) \Subset U'\Subset U$. 

Let us recall some basic facts in pluripotential theory. If $K$ is a compact subset of $U$, for a positive current $\Theta$ of bidegree $(q,q)$, we define the \emph{trace measure}  of $\Theta$  by $\Theta \wedge \omega^{k-q}$. And its mass on $K$ is:
$$   \|\Theta\|_K=\int_K \Theta \wedge \omega^{k-q} . $$ 
We say that a current $S$ satisfies the Chern-Levine-Nirenberg (CLN for short) inequality if for any $(1,1)$ positive closed current $T$ then the positive closed current $T \wedge S$ is well defined and if for any $K \Subset L$ compact subsets, then there exists a constant $C_{K,L,S}>0$ such that:
$$ \|T \wedge  S \|_K \leq C_{K,L,S} \|T\|_L. $$ 
If $S$ is a positive closed current of bidegree $(p,p)$ and $u$ is a bounded psh function, then the current $i\partial \bar{\partial} u \wedge S:=i\partial \bar{\partial}(uS)$ is a well defined positive closed current. Hence, if $u_1,\dots,u_l$ are bounded psh functions, then the current $i\partial \bar{\partial} u_1 \wedge \ldots \wedge i\partial \bar{\partial} u_l \wedge S$ is  a well defined positive closed current. Furthermore, if $K \Subset L$ are compact subsets, then there exists a constant $C_{K,L}>0$ such that the following Chern-Levine-Nirenberg inequality holds:
$$ \|i\partial \bar{\partial} u_1 \wedge \ldots \wedge i\partial \bar{\partial} u_l \wedge  S \|_K\leq C_{K,L}\|u_1\|_{L^\infty} \dots \|u_l\|_{L^\infty} \|S\|_L. $$   
In particular, if $S$ satisfies the CLN inequality, so does $i\partial \bar{\partial} u_1 \wedge \ldots \wedge i\partial \bar{\partial} u_l \wedge S$. Finally, recall that if $(u_n)$ is a uniformly bounded sequence of smooth psh functions decreasing to $u$, then $(i\partial\bar{\partial}u_n\wedge S)$ converges to $i\partial\bar{\partial}u\wedge S$ in the sense of currents. We refer the reader to Demailly's book on that topic (\cite{dem2} pp. 166-172). We have the following lemma (similar to theorem 1 in \cite{Oka1}):
\begin{lemme}\label{PoincaréSobolev}
Let $u$ be a bounded psh function on $U$ and $S$ a positive closed current of bidegree $(p,p)$. Then for $f\in H_S$ we have that $f$ is in $L^2(S\wedge \omega^{k-p-1}\wedge i\partial\bar{\partial}u)$ with:
$$\|f\|^2_{u,S}:=\int_U f^2 S\wedge \omega^{k-p-1} \wedge i\partial\bar{\partial}u \leq 8\|u\|_\infty \|f\|_S^2.$$
\end{lemme}
\emph{Proof} Assume first that $u$ is smooth and $f$ is smooth with compact support in $U$. Denote by $\widetilde{S}$ the current $S\wedge \omega^{k-p-1}$. We consider the norm $|||f|||^2=\int_U f^2 i\partial u\wedge \bar{\partial}u \wedge \widetilde{S}$. Then by Stokes formula:
$$ |||f|||^2 =-2\int_U if\partial f \wedge  u\bar{\partial}u \wedge \widetilde{S} -\int_U f^2 u i\partial\bar{\partial}u\wedge \widetilde{S}.$$
Using the Cauchy-Schwarz inequality:
$$|||f|||^2 \leq 2\|u\|_\infty |||f||| \  \|f\|_S+\|u\|_\infty \|f\|^2_{u,S}.$$
On the other hand:
$$\|f\|^2_{u,S}=-2\int_U i f\partial f \wedge \bar{\partial}u \wedge \widetilde{S} \leq 2|||f||| \  \|f\|_S.$$
Using the two last inequalities, we get:
$$ |||f|||^2 \leq 2\|u\|_\infty |||f||| \ \|f\|_S+2\|u\|_\infty |||f||| \  \|f\|_S.$$
So, we have $|||f|||\leq 4 \|u\|_\infty \|f\|_S$, and so:
$$\|f\|^2_{u,S} \leq 8 \|u\|_\infty \|f\|_S,$$
which is what we want. The general case follows by approximation (first of $u$ by a decreasing sequence of smooth psh functions then of $f$). $\Box$\\

We will need the following definition:
\begin{defi}
Let $S$ be a positive closed current of bidegree $(p-1,p-1)$. We say that $S$ is \emph{$W^*$-regular} if $S$ satisfies the CLN inequality and if the canonical map $\varphi\mapsto \chi \varphi$ from $W^*$ to $H_S$ which is defined for $\varphi$ smooth can be extended to $W^*$ as a bounded linear map which is continuous for the weak topology on $W^*(U')$ and for the weak Hilbert space topology on $H_S$.  
\end{defi}
Recall that a function is continuous for the weak topology if the image of a weakly converging sequence is weakly converging. This notion is interesting here because the weak topology we consider on $W^*$ is not the usual Banach space one (in the case of the weak Banach space topology, being weakly continuous and being strongly continuous are equivalent).

 By lemma \ref{regularization}, any function $\varphi$ in $W^*$ is a limit of a sequence of smooth functions which converges weakly in $W^*(U')$ so the extension is unique. Furthermore, provided that the map $\varphi\mapsto \chi \varphi$ is bounded, a diagonal extraction shows that if the image of any weakly converging sequence of smooth functions is weakly converging then the image of any weakly converging sequence in $W^*$ is weakly converging and thus $S$ is $W^*$-regular. The notion of $W^*$-regularity is interesting because of the following lemma:
\begin{lemme}
Let $S$ be a $W^*$-regular positive closed current of bidegree $(p,p)$ for $p\leq k-2$. Then $S \wedge i\partial \bar{\partial}u$ is $W^*$-regular if $u$ is a bounded psh function. 
\end{lemme}
\emph{Proof.} Denote $S\wedge i\partial \bar{\partial}u$ by $\widetilde{S}$. We know it satisfies the CLN inequality.

We first check that the canonical application $\varphi \to \chi \varphi$ from $W^*$ to $H_{\widetilde{S}}$ is bounded for smooth functions. So, let $\varphi$ be a smooth function in $W^*$ and $T_\varphi$ a positive closed current such that $i\partial \varphi \wedge \bar{\partial}\varphi\leq T_\varphi$.  Then, we have:
\begin{eqnarray*}
\|(\chi \varphi)\|_{\widetilde{S}}^2 \leq 2\int_U \chi^2 i\partial  \varphi\wedge \bar{\partial} \varphi \wedge \widetilde{S}\wedge \omega^{k-p-2}+ 2\int_U \varphi^2 i\partial \chi\wedge \bar{\partial} \chi \wedge \widetilde{S}\wedge \omega^{k-p-2}
\end{eqnarray*}
The first integral of the right hand side is less than $2\int_U  \chi^2 T_\varphi \wedge \widetilde{S}\wedge \omega^{k-p-2}$ which is bounded by the CLN inequality. For the second integral, observe that, choosing some non negative smooth function $\chi_1$ with compact support on $U$ and equal to  1 on the support of $\chi$, we have that:
$$\int_U \varphi^2 i\partial \chi\wedge \bar{\partial} \chi \wedge \widetilde{S}\wedge \omega^{k-p-2} \leq C\int_U (\chi_1 \varphi)^2 i\partial\bar{\partial}u \wedge S\wedge \omega^{k-p-2},$$
 for $C$ large enough depending on $\chi$. Since $S$ is $W^*$-regular, so is $S\wedge\omega^{k-2} $ and we can apply lemma \ref{PoincaréSobolev} to $f=\chi_1 \varphi$. Combining the two estimates, we get:
 $$ \|(\chi \varphi)\|_{\widetilde{S}} \leq A_1 \|\varphi \|_{W^*(U')},$$
 where $A_1$ is a constant that depends only on $S$ and $u$.  \\

Now, we prove the continuity. Let $f$ be a smooth form with compact support in $U$ and $(\varphi_n)$ a sequence in $W^*$ weakly converging in $W^*(U')$. Smooth functions with compact support in $U$ are dense in $H_{\widetilde{S}}$ by definition. So in order to show that $(\chi \varphi_n)$ is weakly converging in $H_{\widetilde{S}}$, it is enough to check that $(\langle \chi\varphi_n, f \rangle_{\widetilde{S}}) $ converges to a value $\langle g, f\rangle_{\widetilde{S}}$ where $g\in H_{\widetilde{S}}$ does not depend on $f$. 

Since $f$ is smooth, there is a $C>0$ such that $i\partial\bar{\partial}f+C\omega\geq 0$. Define $S_1:=S\wedge (i\partial \bar{\partial}f+C\omega)$ and $S_2=S\wedge C\omega$, it is clear that $S_1$ and $S_2$ are $W^*$-regular. Choose some non negative smooth function $\chi_1$ with compact support on $U'$ and equal to 1 on the support of $\chi$. Since $u$ is a bounded psh function, it belongs to $W^*$, so $\chi_1 u$ defines an element of $H_{S'}$ for $S'$ $W^*$-regular. We claim that:
$$\langle \chi \varphi_n, f \rangle_{\widetilde{S}}=  \langle \chi \varphi_n, \chi_1 u \rangle_{S_1}-\langle \chi \varphi_n, \chi_1 u \rangle_{S_2}.  $$
 Indeed, if $u$ is smooth, it is clear by Stokes formula and the general case follows by density since $S_1$ and $S_2$ are $W^*$-regular. The right-hand side shows that $(\langle \chi \varphi_n, f \rangle_{\widetilde{S}})$ converges to the well defined value $\langle \chi \varphi, \chi_1 u \rangle_{S_1}-\langle \chi \varphi, \chi_1 u \rangle_{S_2}$ (indeed, $S_1$ and $S_2$ are $W^*$-regular). By Cauchy-Schwarz inequality and the above bound on the norm of $\chi \varphi_n$, we have that $| \langle \chi \varphi_n, f \rangle_{\widetilde{S}}| \leq A_2 \| \varphi_n\|_{W^*(U')} \|f\|_{\widetilde{S}}$ (where $A_2=\sqrt{A_1}>0$ is a constant that does not depend on $(\varphi_n)$ and $f$). So the mapping $f\mapsto  \lim_{n \to \infty} i\langle \chi \varphi_n, f \rangle_{\widetilde{S}}$ defines a continuous linear form on $H_{\widetilde{S}}$. By Riesz theorem, $(\chi \varphi_n)$ converges weakly in $H_{\widetilde{S}}$ to an element that does not depend on the choice of the sequence $(\varphi_n)$. We still denote by $\chi \varphi$ that element and we have the bound $\|\chi \varphi \|_{\widetilde{S}}\leq A_2 \lim_{n\to \infty} \| \varphi_n\|_{W^*(U')}$. Finally, choosing for $(\varphi_n)$ the sequence in lemma \ref{regularization}, we get that $\|\chi \varphi \|_{\widetilde{S}} \leq A_2 \| \varphi\|_{W^*(U)}$. That completes the proof. $\Box$\\
 
In particular, by induction, we get that any $S=i\partial \bar{\partial} u_1\wedge \ldots \wedge i\partial \bar{\partial} u_{k-1}$ with $u_1,\dots, u_{k-1}$ bounded psh functions on $U$ is $W^*$-regular. And the above proof show that there exists a constant $A$ that only depends on the $L^\infty$ norm of each $u_l$ such that $\| \chi \varphi \|_S \leq A \| \varphi \|_*$. In uniformly convex spaces and thus in Hilbert spaces, there is the classical theorem that will be of use:
\begin{theorem}
Let $A$ be a uniformly convex Banach space and let $(a_n)$ be a sequence in $A$ weakly converging to $a$. Then there is a subsequence $(a_{n_l})_l$, such that the sequence $(\frac{1}{j}\sum_{l=1}^{j} a_{n_l})_j$ converges strongly to $a$. 
\end{theorem}
In particular, for a $W^*$-regular current $S$, we can find a sequence $(\varphi_n)$ converging weakly in $W^*$ and \emph{strongly} in $H_S$ to $\varphi$, of course this sequence depends on $S$ in general. We can assume that this sequence is smooth. Finally, we have proved the following theorem:
\begin{theorem}\label{key}
Let $S=(i\partial \bar{\partial} u)^{k-1}$ where $u$ is a bounded psh function on $U$. Then $S$ is $W^*$-regular. Consequently for $\varphi$ in $W^*$, then $\chi \varphi$ is in $L^2((i\partial \bar{\partial} u)^{k})$ and there is a constant $A$ that depends only on the $L^\infty$ norm of $u$ such that $\| \chi \varphi \|_S \leq A \| \varphi \|_*$. 

Finally, for any sequence  $(\varphi_n)$ converging weakly in $W^*$, there is a subsequence $(\varphi_{n_l})_l$ such that the Cesàro mean $(\chi \frac{1}{j}\sum_{l=1}^{j} \varphi_{n_l})_j$ converges strongly in $H_S$.
\end{theorem} 
Following the results of the appendix, we can now define the functional capacity ${\rm cap}_S$ for $S$ in theorem $\ref{key}$. In particular, the element of $H_S$ are well-defined up to a set of $S$-capacity zero. In particular, for $\varphi$ in $W^*$, the element $\chi \varphi$ of $H_S$ admits a quasi-continuous representative for the capacity ${\rm cap}_S$. This will be useful since we now by \cite{FO1} and \cite{FO2} that if a set is of $S$-capacity equal to zero for all $S$ above, then it is pluripolar. A difficulty is that the quasi-continuous representative depends a priori on $S$ since the converging sequence depends on $S$. 
\begin{Remark} \rm The distribution $i\partial \bar{\partial} \varphi \wedge S := f\mapsto -\langle f, \varphi \rangle_S$ on $V\subset supp(\chi)$ is well defined and of order 1. By Stokes formula, it coincides with the usual definition if $\varphi$ is smooth. Furthermore, if $\varphi_n\rightharpoondown \varphi$ then $i\partial \bar{\partial} \varphi_n \wedge S \to i\partial \bar{\partial} \varphi \wedge S$.
\end{Remark}
 
\subsection{Quasi-continuity} 
Now we want to prove that the functions in $W^*$ are quasi-continuous for the Bedford-Taylor capacity. Recall that $U$ is strongly pseudoconvex so $U= \big\{ \psi <0 \big \}$ where $\psi$ is a smooth strictly psh function on a neighborhood of $\overline{U}$. We will use the intermediate space $W_\infty^*$  consisting of function $\varphi$ in $W^*$ such that there is a current $T_\varphi$ with bounded potential such that $i\partial \varphi \wedge \bar{\partial} \varphi\leq T_\varphi$. Recall that a psh function $v$ is a \emph{potential} of a positive closed current $T$ of bidegree (1,1) if $i\partial\bar{\partial} v=T$. We put on  $W_\infty^*$ the norm $\| \varphi \|^2:= \|\varphi\|_{L^2}^2+ \inf\{\| v_\varphi\|_{L^\infty}| \ v_\varphi \ \text{potential of} \ T_\varphi \}$. All the results of section 1 are still true for $W_\infty^*$: it is a Banach space, and the weak convergence defined as for $W^*$ enjoys the same compactness property. 
\begin{lemme}\label{regularization2}
 Let $\varphi \in W_\infty^*\cap L^\infty$.Then $\chi \varphi$ is in $\in W_\infty^*\cap L^\infty$ and if $\varphi_n \rightharpoondown \varphi$ in $W_\infty^*$ then $\chi\varphi_n \rightharpoondown \chi\varphi$ in $W_\infty^*$. 
 
 Furthermore, there exist sequences $(\varphi_n)$ of smooth functions with $\chi\varphi_n \rightharpoondown \chi\varphi$ with $T$ and $T_n$ satisfying $ i \partial (\chi\varphi) \wedge \bar{\partial} (\chi\varphi) \leq T$ and $i \partial (\chi\varphi_n) \wedge \bar{\partial} (\chi\varphi_n) \leq T_n $ such that the bounded potentials $v_n$ for $T_n$ are decreasing to the bounded potential $v$ for $T$.
\end{lemme}\label{3.6}
\emph{Proof.} Let $\varphi$ as above, and $T_\varphi$ be such $ i \partial\varphi \wedge \bar{\partial} \varphi \leq T_\varphi$. We have the bounds: 
\begin{eqnarray*}
i \partial (\chi\varphi) \wedge \bar{\partial} (\chi\varphi) &\leq& 2i\chi^2 \partial \varphi \wedge \bar{\partial} \varphi+2i|\varphi|^2 \partial \chi\wedge \bar{\partial}\chi\\
																														 &\leq& C_1T_\varphi+C_2 \|\varphi\|_{L^\infty} \omega,						
\end{eqnarray*}
where $C_1$ and $C_2$ are positive constants that depend only on $\chi$. That gives the first part of the lemma.

Now we take $T= C_1T_\varphi+C_2 \|\varphi\|_{L^\infty} \omega$ and $(\varphi_n)$ as in lemma \ref{regularization} which gives the second part part of the lemma. $\Box$ \\

\noindent So for $\varphi$ in $W_\infty^*\cap L^\infty$, we consider  a sequence $(\varphi_n)$ as in the above lemma. Observe that taking a subsequence or a Cesàro mean do not change the fact that the potentials are decreasing (that is simply because a subsequence and a Cesàro mean of a decreasing sequence are still decreasing). We let $u$ be a psh function on $U$ with $0\leq u \leq 1$, and $\varphi_n\rightharpoondown \varphi$ as in the above lemma. We will need the following lemma:
\begin{lemme}
With the notations of lemma \ref{3.6}, for all $j\leq k$, there exists $C>0$ which depends only on $\varphi$ such that:
$$ \int_U |\chi \varphi -\chi \varphi_n|^2 (i \partial\bar{\partial}u)^k \leq C \Big( \int_U |\chi \varphi -\chi \varphi_n|^2 (T+T_n)^j\wedge(i \partial\bar{\partial}u)^{k-j} \Big)^{\frac{1}{2^j}}.$$ 
\end{lemme}
\emph{Proof.} Since we assume that $u$ and the potentials of $T$ and $T_n$ are bounded, the previous integrals make sense by lemma \ref{PoincaréSobolev} ($\xi\varphi-\xi\varphi_n$ is in $L^2((T+T_n)^j\wedge(i \partial\bar{\partial}u)^{k-j})$ for $j\geq 0$). We prove the claim of the lemma by induction on $j$. For $j=0$ there is nothing to prove, so assume the claim holds for $j$. The Stokes formula implies:
\begin{eqnarray*}
 \lefteqn{\int_U |\chi \varphi -\chi \varphi_n|^2 (T+T_n)^j\wedge(i \partial\bar{\partial}u)^{k-j} =} \\
 & & -2 \int_U (\chi \varphi -\chi \varphi_n)i\partial(\chi \varphi -\chi \varphi_n)\wedge\bar{\partial}u\wedge (T+T_n)^j\wedge(i \partial\bar{\partial}u)^{k-j-1},
\end{eqnarray*}
and by Cauchy-Schwarz inequality and the fact that $i\partial(\chi \varphi -\chi \varphi_n)\wedge \bar{\partial}(\chi \varphi -\chi \varphi_n)\leq 2(T+T_n)$, it is bounded by:
\begin{eqnarray*}
 2\Big(\int_U |\chi \varphi -\chi \varphi_n|^2 2(T+T_n)^{j+1}\wedge(i \partial\bar{\partial}u)^{k-j-1}\Big)^{\frac{1}{2}} \times \\
 \Big(\int_{Supp(\chi)} i\partial u\wedge\bar{\partial}u \wedge(T+T_n)^{j}\wedge(i \partial\bar{\partial}u)^{k-j-1}\Big)^{\frac{1}{2}}.
 \end{eqnarray*}
The last term of the product is less than to $2^{-1}\int_{Supp(\chi)}  i\partial\bar{\partial}(u^2) \wedge(T+T_n)^{j}\wedge(i \partial\bar{\partial}u)^{k-j-1}$ since $i\partial u \wedge \bar{\partial}u \leq i2^{-1}\partial\bar{\partial}(u^2)$ since $u\geq 0$ (expand the right-hand side of the inequality). And that quantity is bounded by a constant independent of $u$ by the CLN inequality. That prove the claim for $j+1$. $\Box$\\

We now prove the following lemma which is the key point for the proof of the quasi-continuity.
\begin{lemme}\label{convergence}
Let  $\varphi \in W_\infty^*\cap L^\infty$. Then there exists a sequence of smooth functions $\varphi_n \rightharpoondown \varphi$ in $W_\infty^*\cap L^\infty$ such that for any psh function $0\leq u \leq 1$ on $U$, we have:
 $$ \int_U |\chi \varphi -\chi \varphi_n|^2 (i \partial\bar{\partial}u)^k \to 0,$$
where the convergence is uniform in $u$.
\end{lemme} 
\emph{Proof.} We apply the previous lemma for $j=k$ and for $(\varphi_n)$ in lemma \ref{regularization2}. Proving that $\int_U  |\chi \varphi -\chi \varphi_n|^2(T+T_n)^{k}$ goes to zero will give the lemma. We expand in the integral and we obtain terms of the form $\int_U |\chi \varphi -\chi \varphi_n|^2T^{k-j}\wedge T_n^j$. We prove by induction on $j$ that we can find a sequence $\varphi_n \rightharpoondown \varphi$ such that the integrals $\int_U |\chi \varphi -\chi \varphi_n|^2T^{k-l}\wedge T_n^l$ goes to zero for $l\leq j$. Theorem \ref{key} applied to $S=T^{k-1}$ gives for $j=0$, that possibly after taking a subsequence and a Cesàro mean, the integral $\int_U |\chi \varphi -\chi \varphi_n|^2T^k$ goes to zero. Assume the claim hold for $j-1$, so we have a sequence $\varphi_n \rightharpoondown \varphi$ such that the integrals $\int_U |\chi \varphi -\chi \varphi_n|^2T^{k-l}\wedge T_n^l$ goes to zero for $l\leq j-1$. We write $T_n=T+T_n-T$ in $\int_U |\chi \varphi -\chi \varphi_n|^2T^{k-j}\wedge T_n^j$, so we get the sum:
$$\int_U |\chi \varphi -\chi \varphi_n|^2T^{k-j}\wedge (T_n-T)\wedge T_n^{j-1}+ \int_U |\chi \varphi -\chi \varphi_n|^2T^{k-j+1}\wedge T_n^{j-1}.$$
The second term of the right-hand side goes to zero by the claim for $j-1$. Recall that we call $v$ and $v_n$ the potentials of $T$ and $T_n$. Observe that by Stokes formula:
\begin{eqnarray*}
\lefteqn{\int_U |\chi \varphi -\chi \varphi_n|^2T^{k-j}\wedge (T_n-T)\wedge T_n^{j-1}=}\\
& & -2\int_U (\chi \varphi -\chi \varphi_n)i\partial(\chi \varphi -\chi \varphi_n) \wedge\bar{\partial}(v-v_n) \wedge T^{k-j}\wedge T_n^{j-1}.
\end{eqnarray*}
Once again, using Cauchy-Schwarz inequality, we can bound the last term by:
$$ 2\Big(\int_U |\chi \varphi -\chi \varphi_n|^2 2(T_n+T)\wedge T^{k-j} \wedge T_n^{j-1}\Big)^{\frac{1}{2}}\Big(\int_U i\partial(v_n-v)\wedge \bar{\partial}(v_n-v) \wedge T^{k-j} \wedge T_n^{j-1}\Big)^{\frac{1}{2}}.$$
The first term is bounded by lemma \ref{key}. For the second one, observe that by a standard argument of max construction, we can assume that there exists a constant $B$ depending only on the $L^\infty$ norm of the psh function $v$ such that the functions $v$ and $v_n$ coincide with $B \psi$ in a neighborhood $W$ of $\partial U$ that does not depend on $n$ (e.g. \cite{dem2} p. 170). So, by Stokes formula, the second term is equal to: 
$$ -\int_U (v_n-v) (T_n-T)\wedge T^{k-j} \wedge T_n^{j-1}. $$
That term is equal to $ -\int_U  \chi_1 (v_n-v) (T_n-T)\wedge T^{k-j} \wedge T_n^{j-1}$ for $\chi_1$ a non negative smooth function with compact support on $U$ and equal to 1 on $U\backslash W $. And since the $v_n$ are uniformly bounded and decreasing to $v$, we know that $(v_n-v) (T_n-T)\wedge T^{k-j} \wedge T_n^{j-1}$ converges to $0$ in the sense of currents by theorem 3.7 in \cite{dem2} p-170. That gives the claim. The result follows for $j=k$. $\Box$\\

We want to define functions in $W^*$ up to a pluripolar set, that is we want to find a representative of $\varphi$ in $W^*$ defined up to a pluripolar set. Such a definition would not be of much interest if two such representatives differ on a set bigger than pluripolar. In order to prove that the definition of the representative is meaningful, we will need the notion of pluri-fine topology.
Recall that the \emph{pluri-fine topology} is the coarsest topology for which the psh functions are continuous (pluri-finely and pluri-fine will refer to this topology). This topology is strictly finer than the usual topology. By definition, complete pluripolar sets are pluri-finely closed.
For a Borel set $E$, we have that ${\rm cap}_{BT}(E)= \int_U (dd^c v_E)^k$ where $v_E$ is the \emph{extremal function} associated to $E$. It is defined by $v_E= (\sup \big\{u , \ u \ \text{psh on} \ U, \ u=-1 \ \text{on} \ E, \ u\leq 0 \ \text{on} \ U \big\})^*$  where $f^*$ denotes the upper semi continuous regularization of a function $f$. The function $v_E$ is psh, non positive and equal to $-1$ on $E$ apart from a pluripolar set \cite{BT1}. 
Let $E\subset U$ be a Borel set, and denote by $\overline{E}^f$ its pluri-fine closure. Then ${\rm cap}_{BT}(E)={\rm cap}_{BT}(\overline{E}^f)$ because the extremal function of $E$ and $\overline{E}^f$ are equal by definition.
\begin{proposition}
Let $f$ be a function which is quasi-continuous for the Bedford-Taylor capacity. Then it is pluri-finely continuous outside a pluripolar set.
\end{proposition}
\emph{Proof.} Take $f$ as above. By definition, for all $n\geq1$, there exists an open set $V_n$ in $U$ of capacity less than $n^{-1}$ such that $f$ is quasi-continuous on $U\backslash V_n$. Considering $\cap_{j\leq n} V_j$, we can assume that the sequence $(V_n)$ is decreasing. By restriction, $f$ is continuous on the pluri-open set $U\backslash\overline{V_n}^f$ and we know from above that ${\rm cap}_{BT}(\overline{V_n}^f) \leq n^{-1}$. Consider the set $P:=\cap_n \overline{V_n}^f$, it is pluripolar because its capacity is equal to zero (it is less than $n^{-1}$ for all $n$). Then, the function $f$ is pluri-finely continuous on its complement because for $x \notin P$, then $x$ is in the pluri-finely open set $U\backslash\overline{V_n}^f$ for some $n$ and it is continuous there. $\Box$
\begin{proposition}\label{unicity}
Let $f$ be a function pluri-finely continuous outside a pluripolar set and vanishing almost everywhere, then it vanishes outside a pluripolar set.
In particular, two quasi-continuous representatives of a function are equal outside a pluripolar set. 
\end{proposition}
\emph{Proof.} Let $f$ be as above. Observe that pluri-finely open sets are either of positive Lebesgue measure or empty. Indeed, we only have to check that $u^{-1}(\{x>c\})$ is of positive Lebesque measure or empty for $u$ psh and $c \in \mathbb{R}$. That is the case by the mean value inequality. In particular, the pluri-finely open set $\{f\neq 0\}$ is empty. The rest of the proposition follows. $\Box$ \\

We can now prove the quasi-continuity result.
\begin{theorem}\label{quasi-continuity2}
Let $\varphi\in W^*$ then there exists a representative of $\varphi$ which is quasi-continuous for the Bedford-Taylor. In particular, this representative is pluri-finely continuous outside a pluripolar set and two such representatives coincide outside a pluripolar set.
\end{theorem}
\emph{Proof.} 
 Consider an increasing sequence of open sets $V_i\Subset U$ with $\cup_i V_i= U$. Assume that for all $i$ we can find a representative $\widetilde{\varphi}^i$ of $\varphi$ that is quasi-continuous on $V_i$. Since $\widetilde{\varphi}^j$ is quasi-continuous on $V_i$ for $j\geq i$, we can remove a countable union of pluripolar sets on which the different representatives do not coincide. That way, we can take a representative $\widetilde{\varphi}$ which is quasi-continuous on all $V_i$. Then for $\varepsilon>0$, $\widetilde{\varphi}$ is continuous on  $V_i\backslash G_i$ where ${\rm cap}_{BT}(G_i) \leq \varepsilon 2^{-i}$. So $\widetilde{\varphi}$ is continuous outside $G=\cup G_i$. Since ${\rm cap}_{BT}$ is subadditive, we have ${\rm cap}_{BT}(G)\leq\varepsilon$. Thus $\widetilde{\varphi}$ is quasi-continuous on $U$. So it is sufficient to prove that $\varphi$ admits a representative that is quasi-continuous on a open set $V$ with $V\Subset U$. Choose a smooth function $\chi$ with compact support in $U$ so that $V\subset \{\chi=1\}$.
 
We prove that $\chi \varphi$ is quasi-continuous on $U$. The proof is in three steps. First we assume that $\varphi \in W^*_\infty \cap L^{\infty}$, then we extend the result to functions in $W^* \cap L^{\infty}$ and finally we prove the general case.\\

\noindent {\bf{Step 1.}} By lemma \ref{convergence}, we choose $(\varphi_n)$ smooth, weakly converging to $\varphi \in  W^*_\infty \cap L^{\infty}$ such that $\int_U |\chi \varphi -\chi \varphi_n|^2 (i \partial\bar{\partial}u)^k$ goes to zero uniformly in $u$. For $\alpha>0$, the sets $E_n:=\big\{ |\varphi_n -\varphi| \geq \alpha \big \}$ satisfies:
$${\rm cap}_{BT}(E_n) \leq \sup \frac{1}{\pi^k} \Big\{\frac{1}{\alpha} \int_U |\chi \varphi -\chi \varphi_n|^2 (i \partial\bar{\partial}u)^k , \ u \ \text{psh}, \ 0\leq u\leq 1 \Big\},$$
 since $\pi dd^c=  i\partial \bar{\partial}$. Hence the sequence $({\rm cap}_{BT}(E_n))$ goes to zero. Taking a subsequence, we can assume that the sequence $(\varphi_n)$ satisfies ${\rm cap}_{BT}(\{ |\varphi_n-\varphi|>2^{-n-1} \})< 2^{-n-1}$. So, $\sum_n (\chi\varphi_{n+1}-\chi\varphi_{n})+\chi\varphi_1$ converges uniformly outside the open set $\cup_{n\geq j} \{ |\varphi_n-\varphi_{n+1}|>2^{-n} \}$, which is of BT-capacity less than $2^{1-j}$ (recall that the BT-capacity is subadditive). That gives the first step.\\

\noindent {\bf{Step 2.}} Now, let $\varphi \in  W^* \cap L^{\infty}$ so there exists a psh function $v<0$ such that $i\partial \varphi \wedge \bar{\partial}\varphi \leq i  \partial\bar{\partial}v$. In order to simplify the notations, assume $\|\varphi\|_{L^\infty} \leq 1$. For $N>0$, let $v_N:=\sup(v,-N)$ and $\varphi_N=N^{-1}(N+v_N) \varphi$ (that way, $\varphi_N$ is equal to zero where $v<-N$).
We want to show that $\varphi_N$ is in $ W^*_\infty \cap L^{\infty}$. Assume first that $\varphi$ and $v$ are smooth. We have the bound:
$$i\partial \varphi_N \wedge \bar{\partial}\varphi_N \leq 2i\Big((\frac{N+v_N}{N})^2\partial \varphi \wedge \bar{\partial}\varphi+ |\varphi|^2\partial (\frac{N+v_N}{N}) \wedge \bar{\partial}(\frac{N+v_N}{N})\Big).$$
 By definition of $\varphi$, we have $(1+v_N/N)^2 i\partial \varphi \wedge \bar{\partial}\varphi  \leq (1+v_N/N)^2 i\partial\bar{\partial}v$ and we have $(1+v_N/N)^2 i\partial\bar{\partial}v \leq i\partial\bar{\partial}v_N$. Indeed, it is immediate on the open sets $\{v>-N\}$ and the interior of $\{v \leq -N\}$ and since the left-hand side is a continuous current, it does not give mass to $\partial \{v \leq -N\}$ so the inequality holds in the sense of currents. It implies that:
$$i\partial \varphi_N \wedge \bar{\partial}\varphi_N\leq 2i\partial\bar{\partial}v_N+\frac{2}{N^2} i\partial v_N \wedge \bar{\partial}v_N.$$
 Now $2i\partial v_N \wedge \bar{\partial}v_N \leq i\partial\bar{\partial}(v_N+N)^2$, so $\varphi_N$ is in $W_{\infty}^* \cap L^{\infty}$ with a control of the norm depending only on $N$. More precisely, we have the bound $\|\varphi_N\|_{L^\infty}\leq 1$ and the current satisfying (\ref{definition}) for $\varphi_N$ has a potential taking values in $[-2N,2]$. Taking a weak limit as in lemma \ref{regularization2}, we deduce that $\varphi_N$ is always in $W^*_\infty \cap L^\infty$ for $\varphi$ in $W^* \cap L^\infty$. \\
That defines a sequence $(\varphi_N)$ in $W^*_\infty \cap L^{\infty}$ which converges weakly to $\varphi$ in $W^*\cap L^{\infty}$. Let $E_\infty=\{ v=-\infty\}$ and for each $n\in \mathbb{N}$, let $P_n$ be the set of points where a fixed quasi-continuous representant $\widetilde{\varphi}_n$ of $\varphi_n$ is not defined. Let $P:=E_\infty \cup(\cup_n P_n)$ so $P$ is pluripolar. Let $x\notin P$. Then for $N$ large enough, we have that $x \in \{ v > -N\}$, so extend the definition (which we know stands a priori almost everywhere) of $\varphi$ at $x$ by:
$$\widetilde{\varphi}(x):= \frac{N}{N+v_N(x)}\widetilde{\varphi}_N(x).$$
It is crucial here that $\widetilde{\varphi}(x)$ does not depend on the choice of $N$. That is the case because if $N\geq N'$ then the function:  $$g_n:=\Big(\frac{N}{N+v_N}\widetilde{\varphi}_N-\frac{N'}{N'+v_N'}\widetilde{\varphi}_{N'}\Big)\frac{N'+v_N'}{N'}$$
is well defined and equal to: 
$$ \widetilde{\varphi}_N(x)\frac{N(N'+v_N')}{N'(N+v_N)}- \widetilde{\varphi}_{N'}$$
which is pluri-finely continuous outside a pluripolar set as product and sum of functions pluri-finely continuous there. It vanishes almost everywhere so it vanishes outside a pluripolar set by proposition \ref{unicity}. So removing a countable union of pluripolar sets if necessary, we can define $\varphi$ quasi-everywhere. 

Now, for $\varepsilon >0$, take $N$ large enough so that ${\rm cap}_{BT} \{v<-N \}\leq \varepsilon$ (see \cite{BT1}). On $\{v\geq -N \}$ $\chi\widetilde{\varphi}$ is given by a function which is quasi-continuous. That gives the result for step 2.\\

\noindent {\bf{Step 3.}} Let now be $\varphi \in W^*$, it is sufficient to assume that $\varphi \geq 0$ since we can write $\varphi= \varphi^+-\varphi^-$. Let $N\in \mathbb{N}$, we define $\varphi_N$ by $\varphi_N(x)= \inf (\varphi,N)$ so $\varphi_N \in W^*\cap L^\infty$ . By Step 2, we know that $\varphi_N$ is quasi-continuous, hence it is defined up to a pluripolar set (we make the identification between $\varphi_N$ and one of its quasi-continuous representative). Let $j>0$, we remark that $\varphi_N(x)=\inf(\varphi_{N+j}(x),N)$ for almost every $x$, so it is true outside a pluripolar set by proposition \ref{unicity}. As before, let $P$ be the pluripolar set consisting of all the points where the different functions $\varphi_N$ are not well defined and where the equalities $\varphi_N(x)=\inf(\varphi_{N+j}(x),N)$ do not hold. Let $x\notin P$, we have that $(\varphi_N(x))_N$ is an increasing sequence constant for $N$ large enough, so we defined $\varphi(x)=\lim_N \varphi_N(x)$.  Let $F_N:=(\{\varphi \geq N\}\cup P)\cap supp(\chi)$. Let $\chi_1$ be a smooth function with compact support in $U$ such that $supp(\chi) \Subset \{\chi_1=1 \}$. Then ${\rm cap}_{BT}(F_N)\to 0$ since for $0\leq u\leq 1$ a psh function: 
$$\int_{F_N} (i\partial\bar{\partial} u)^k\leq \int_U  (\frac{\chi_1 \varphi}{N})^2(i\partial\bar{\partial} u)^k,$$
and we can conclude by lemma \ref{PoincaréSobolev} and theorem \ref{key}.

So $\chi\varphi$ is quasi-continuous on $U\backslash F_N$. Indeed, for $\varepsilon>0$, take $N$ so that ${\rm cap}_{BT}(F_N)\leq \varepsilon$ and take an open set $U_N$ in $V$ such that $F_N\subset U_N$ with ${\rm cap}_{BT}(U_n)\leq 2 \varepsilon$ (this is possible because the Bedford-Taylor capacity is outer regular). Outside $U_N$, $\chi\varphi=\chi\varphi_N$ which is quasi-continuous. That concludes the proof. $\Box$ \\

Since we can take representatives defined up to a pluripolar set, from now on, $\varphi$ will denote a quasi-continuous representative of $\varphi$ for the BT-capacity. Recall that a function in $H_S$ also admits a quasi-continuous representative for the capacity ${\rm cap}_S$ and if $E$ satisfy ${\rm cap}_{BT}(E)=0$ then ${\rm cap}_S(E)=0$ (see \cite{FO1}). So it is natural to ask if $\chi \varphi$ is also a quasi-continuous representative of $\chi \varphi \in H_S$ for all $S$ as in theorem \ref{key}.
\begin{lemme}\label{3.12}
Let $S$ be as in theorem \ref{key}. Let $\varphi$ denote a representative of an element of  $W^*$ quasi-continuous on $U$. Then $\chi \varphi$ is also quasi-continuous for the capacity ${\rm cap}_S$.  
\end{lemme}
\emph{Proof} We follow each step of the previous proof and we check at each step that the quasi-continuous representative for the BT-capacity we defined is also quasi-continuous for the $S$-capacity.

 The lemma holds for $f\in W_\infty^*\cap L^\infty$ by theorem \ref{key}. Indeed, we can take a subsequence and a Cesàro mean of the sequence $(\varphi_n)$ of smooth functions converging to $\varphi$ so that it is strongly converging to $\varphi$ in $H_S$. That gives the result by lemma \ref{lusin2}.

 Now, for $\varphi \in W^*\cap L^\infty$, by construction of the quasi-continuous representative for the Bedford-Taylor capacity, it is sufficicent to check that ${\rm cap}_S(\{v \leq -N \})\to 0$ for $v$.  This is true since, taking $-\log (-v)$ if necessary we can assume that $v$ is in $W^*$ (see the first example of section \ref{examples}) and then by theorem \ref{key}:
 $$ {\rm cap}_S(\{v \leq -N \}) \leq \int_U  i\partial{(\frac{\chi v}{N})} \wedge \bar{\partial}(\frac{\chi v}{N})\wedge S \leq (\frac{A}{N}\|v\|_*)^2.$$ 
And finally, for $\varphi \in W^*$, it is sufficient to check that ${\rm cap}_S(\{ \chi \varphi \geq N\}) \to 0$ which is clear by theorem \ref{key} since:
 $$ {\rm cap}_S(\{ \chi \varphi \geq N\})\leq \int_U  i\partial{(\frac{\chi\varphi}{N})} \wedge\bar{\partial}(\frac{\chi\varphi}{N})\wedge S \leq (\frac{A}{N}\|\varphi\|_*)^2.$$
 That concludes the proof. $\Box$\\
 
\noindent We can now prove the following result on pointwise convergence which will be useful in proving continuity result for the capacity. It can be seen as a weak version of lemma \ref{lusin2}:
\begin{lemme}\label{pointwise}
Let $(\varphi_n)_n$ be a sequence in $W^*$ weakly converging to $\varphi$ in $W^*$. Let $a\in \mathbb{R}$ and let $A$ be a Borel set such that each $\varphi_n$ is equal to $a$ on $A\backslash H_n$ where $H_n$ is a pluripolar set. Then $\varphi$ is equal to $a$ on $A\backslash H$ where $H$ is a pluripolar set.
\end{lemme} 
\emph{Proof.} The result is local so we only prove it on $V\Subset U$ with $\chi$ as above.  We choose $S$ as in theorem \ref{key}. Extracting and using a Cesàro mean, we can assume that the sequence $(\chi\varphi_n)$ is strongly converging on $H_S$. So we can extract a sequence converging outside a set of $S$-capacity zero by lemma \ref{lusin2}. Thus, by lemma \ref{3.12}, the result is true on $A$ minus a set of $S$-capacity equal to zero for all $S$. By \cite{FO1}, we know that such sets are exactly pluripolar sets. $\Box$ \\

Now, for $\varphi \in W^*$, we want to estimate the size of its \emph{Lebesgue set}. Recall that the Lebesgue set of a function is defined as the set of points where the mean value is converging. Let $W$ be a pluri-fine neighborhood of $x \in V$ and $\psi$ a psh function in $V$ with $\psi(x)\neq -\infty$, then it is an easy exercise that: $$\frac{1}{r^{2k}}\int_{ B(x,r) \backslash W} \psi d\lambda \to 0,$$
when $r$ goes to zero (e.g. \cite{Bre} p.79). In particular, a function $\varphi \in W^*$ which is pluri-finely continuous at $x$ such that there are psh functions $\psi_1$ and $\psi_2$ with $\psi_1 \leq \varphi \leq -\psi_2$ on $V$ satisfies $m_{B(x,r)}(f) \to f(x)$ for all $x$ such that $\psi_1(x)$ and $\psi_2(x)$ are finite.  So we have the nice corollary:
\begin{corollaire}
If a function $f$ in $W^*$ is bounded, then the complement of its Lebesgue set is pluripolar.
\end{corollaire} 
In order to extend this result to unbounded functions, the only thing left in our approach is to show that every function in $W^*$ can be locally bounded by two psh functions.

\section{Functional capacity}

From now on, $X$ is a compact Kähler manifold, this case was our primary motivation in this work. We will explain in remark \ref{generalisation} how to extend the result in the local case. We deduce from the above section that an element $\varphi$ in $W^*$ admits a quasi-continuous representative for the Bedford-Taylor capacity that we  will still denote by $\varphi$. Recall that a function is said to be quasi-psh (qpsh for short) if it is locally the sum of a psh function and a smooth function, so it satifies $i\partial \bar{\partial} u+C \omega \geq 0$ for some $C$. The Bedford-Taylor capacity can be generalized to compact Kähler manifold by:
$$ {\rm cap}_\omega(E)= \sup \Big\{ \int_E (i\partial \bar{\partial} u + \omega)^k | \ u \ \text{qpsh}, \ i\partial \bar{\partial} u+\omega \geq 0 , \ 0\leq u \leq 1 \Big\}.$$
We will use some of the results of \cite{GZ1} in this section, and we refer the reader to \cite{dem1}, \cite{DS2} and \cite{GZ1} for basics on qpsh functions. In particular, the capacity ${\rm cap}_\omega$ is comparable with ${\rm cap}_{BT}=\sum {\rm cap}_{BT,U_i}(E\cap U_i)$ where $(U_i)_{i\leq M}$ is a finite covering of $X$ by pseudoconvex open sets and ${\rm cap}_{BT,U_i}$ denotes the Bedford-Taylor capacity of $U_i$ (i.e. there is a $A > 0$ with $(1/A) {\rm cap}_{BT}\leq {\rm cap}_\omega \leq A {\rm cap}_{BT}$). An important fact is that the family $\{ u \ \text{qpsh}| \ i\partial \bar{\partial} u+\omega \geq 0 \ \text{and} \ u \leq 0 \}$ is compact for the $L^1$ norm. A set is \emph{globally pluripolar} if it is contained in the set $\{v=-\infty \}$ for some $v$ qpsh on $X$. It turns out that locally pluripolar sets are in fact globally pluripolar \cite{GZ1}, so we will simply speak of pluripolar sets. \\

We will need an equivalent of Alexander capacity (\cite{Ale}, see also \cite{SW1}) for compact Kähler manifolds which was introduced in \cite{DS2} and developped in \cite{GZ1}. For an open set $U$, we consider the function:
$$V_{U,\omega}(x):=\sup \Big\{u(x) \ \text{qpsh}| \ i\partial \bar{\partial} u+\omega \geq 0 \ \text{and} \ u =0 \ \text{on} \ U \Big\}.$$
Then $V_{U,\omega}$ is in fact qpsh, positive, zero on $U$ and satisfies $i\partial \bar{\partial} u+\omega \geq 0$. Then we define the capacity of an open set by:
$$T_\omega(U):=\exp \Big(-\sup_X(V_{U,\omega})\Big).$$
The following estimate was proven in the local case in \cite{AT1} and the same argument gives:
\begin{theorem}\label{comparaison}
There is a $A>0$ such that for any open set $U$ of $X$:
$$\exp \Big[ -\frac{A}{{\rm cap}_\omega(U)} \Big] \leq T_\omega(U) \leq e.\exp \Big[ -\frac{1}{{\rm cap}_\omega(U)^\frac{1}{k}}  \Big].$$
\end{theorem} 
Now, we want to introduce a functional capacity similar to the classical one (see the appendix or \cite{FZ1}). For a Borel set $E$ in $X$, we define the set $L(E)$ of $W^*$ by
$$L(E):= \Big\{ \varphi \in W^*, \  \varphi \leq -1 \ \text{ a.e on some neighborhood of} \ E, \ \varphi \leq 0 \ \text{on} \ X \Big\}.$$
  We define the capacity $C(E)$ of $E$ by:
\begin{eqnarray*}
 C(E):= \inf \Big\{ \| \varphi \|_{*}^2| \ \varphi \in L(E) \Big\} 
\end{eqnarray*}
For a capacity $c$, let $c_*$ be the \emph{inner capacity} associated to $c$ defined by $c_*(E):=\sup \big\{c(K), K \ \text{compact}, \ K\subset E \big\}$ and $c^*$ the \emph{outer capacity} associated to $c$ defined by $c^*(E):=\sup \big\{c(U), \ U\ \text{open}, E\subset U \big\}$. The capacity is said to be \emph{inner regular} if $c=c_*$ and \emph{outer regular} if $c^*=c$, finally a capacity is \emph{regular} if it is both inner and outer regular. By definition, the capacity $C$ is outer regular. Since $\|\max(\varphi,-1)\|_* \leq \| \varphi\|_*$, it is equivalent to take functions equal to $-1$ on some neighborhood of $E$.  
  We have the following properties:
\begin{proposition}\label{firstproperties}
The capacity $C$ satisfies that:
\begin{enumerate}
 \item for $E \subset F \subset X$, $C(E)\leq C(F)$,
 \item for $(E_i)_i$ a sequence of Borel sets in $X$, $C(\cup_i E_1) \leq \sum_i C(E_i)$,  
 \item for any $E$, one has $C(E)\leq 1$.
 \item for $(K_n)$ a decreasing sequence of compact sets, $\lim C(K_n)\to C(\cap K_i)$
\end{enumerate}
 \end{proposition}
\emph{Proof.} The first item and the fourth item are clear.

 For the second one take $(E_i)_i$ a sequence of Borel sets in $X$ and assume the sum $\sum_i C(E_i)$ converges or there is nothing to prove. For each $i$ let $\varphi_i$ be an element in $L(E_i)$ with $C(E_i) \geq \|\varphi_i \|_*^2 -2^{-i}\varepsilon$ and $h_n =\inf_{i\leq n} \varphi_i$. The sequence $h_n$ is decreasing. Recall by formula (\ref{extremum}) in section \ref{examples} that $\|h_n\|^2_* \leq \sum_{i\leq n} \| \varphi \|^2_*$ so we have the bound:
 $$\|h_n\|_*^2 \leq \sum_i C(E_i) +\varepsilon.$$
  Taking a weak limit gives that $h:=\lim h_n$ is in $W^*$ with $\| h \|^2_* \leq \sum_i C(E_i) +\varepsilon$. Since $h \in L(\cup_i E_i)$ and $\varepsilon$ is arbitrary, the assertion follows. 
 
The item 3 is obtained for $\varphi:=1 \in L(E)$.  $\Box$\\

We now state one of our main result that shows that the capacity $C$ characterizes pluripolar sets.

\begin{theorem}\label{pluripolar}
There exists a constant $B>0$ such that for all Borel set $E$, we have ${\rm cap}_\omega(E) \leq B C(E)$. In particular, $C(E)=0$ if and only if $E$ is pluripolar.   
\end{theorem} 
\emph{Proof.} The first assertion is a restating of the results of theorem \ref{key}. We deduce that $C(E)=0$ implies that $E$ is pluripolar because its BT-capacity is zero. 
Now, let $E$ be a pluripolar set in $X$ so there exists a qpsh function in $X$ with $E \subset \varphi^{-1}(-\infty)$. Subtracting a constant if necessary, we can consider $\psi = -\log (- \varphi)$ which is in $W^*$ with the same poles set as $\varphi$ (see the first example of section \ref{examples}). Since $\psi$ is upper semi continuous and non positive, taking $-\psi / N$ for $N$ large enough, gives $C(E)\leq \varepsilon$ for any $\varepsilon >0$. The proposition follows. $\Box$ \\

The previous result is a generalization of the case of Riemann surfaces where polar sets are sets of capacity equal to zero. The following proposition shows that ${\rm cap}_\omega$ and $C$ defines equivalent capacities.
\begin{proposition}\label{quasi-continuity3}
There exists a continuous function $g:\mathbb{R}^+\to \mathbb{R}^+$ with $g(0)=0$ such that $C\leq  g({\rm cap}_\omega)$. So functions in $W^*$ are quasi-continuous for the functional capacity $C$.
\end{proposition}
\emph{Proof.} Consider an open set $U$ such that ${\rm cap}_\omega (U)\leq \varepsilon$ for $\varepsilon$ small. We know that the function $f_U:=(V_{U,\omega}-\max_X V_{U,\omega})/\max_X V_{U,\omega}$ is equal to $-1$ on $U$ with $-1\leq f_U \leq 0$ and $i\partial \bar{\partial}f_U+ \|V_{U,\omega}\|^{-1}_{\infty}\omega \geq 0$ and we have:
\begin{eqnarray*}
i\partial f_U \wedge \bar{\partial}f_U&=&  -f_U i\partial \bar{\partial} f_U+ \frac{1}{2} i\partial \bar{\partial} (f_U^2) \\
                                      &\leq& i\partial \bar{\partial} f_U + \frac{\omega}{\max_X V_{U,\omega}}+ \frac{1}{2} i\partial \bar{\partial} (f_U^2).
\end{eqnarray*}
The right-hand side is a positive closed current of mass $(\max_X V_{U,\omega})^{-1}$. By theorem \ref{comparaison}, it goes to zero with $\varepsilon$. Since the set $\{ u \ \text{qpsh}| \ i\partial \bar{\partial} u+\omega \geq 0 \ \text{and} \ u \leq 0 \}$ is compact for the $L^1$ norm, we also have that $\|f_U\|_2^2$ goes to zero with $\varepsilon$. That gives the proposition since if $(U_n)$ is a sequence of open set with ${\rm cap}_\omega (U_n) \to 0 $ then $C(U_n)\to 0$. $\Box$    \\

\noindent We now want to show the crucial property $C(\cup E_i)=\lim C(E_i)$ for $E_i$ an increasing sequence of Borel sets. So $C$ is a Choquet capacity. This is interesting because Choquet capacities are regular (see \cite{Cho}, theorem 1). For this, we will need an alternative description of the capacity $C$ that uses the fact that the elements of $W^*$ are defined up to a pluripolar set.

\begin{theorem}
For a Borel set $E$, we have that
$$C(E)=\inf \Big\{\|\varphi\|^2_*, \  \varphi \leq 0, \ \varphi \leq -1 \ \text{on} \ E\backslash H_\varphi, \  H_\varphi \ \text{pluripolar} \Big\}.$$ 
In particular, if $E_i$ is an increasing sequence of Borel sets, then:
 $$C(\cup E_i)=\lim C(E_i).$$ 
Thus $C$ is a Choquet capacity.
\end{theorem}
\emph{Proof.} For a Borel subset $E$, we denote by $C'(E)$ the quantity:
$$C'(E):=\inf \Big\{\|\varphi\|^2_*, \  \varphi \leq 0, \ \varphi \leq -1 \ \text{on} \ E\backslash H_\varphi, \  H_\varphi \ \text{pluripolar} \Big\}.$$ 
Clearly, we have $C'(E) \leq C(E)$ so we only need to prove the other inequality. For that, let $\varphi \in W^*$ be a non positive function, less than $-1$ on $E \backslash H$ where $H$ is a pluripolar set, such that $\| \varphi \|_*^2 \leq C'(E)+\varepsilon_1$. We want to slightly modify $ \varphi$ so that it is in $L(E)$. Adding $\varepsilon \psi$, with $\varepsilon$ small, where $\psi$ is a qpsh function equal to $-\infty$ on $H$ and taking $\max( \varphi+\varepsilon \psi,-1)$ we can assume that $H=\varnothing$. Now, let $\varepsilon_2>0$ be such that $\varphi$ is continuous on the complement of some open set $U$ with $C(U)\leq \varepsilon_2$. Let $\varphi_U$ be in $L(U)$ with $\|\varphi_U \|_* \leq 2\varepsilon_2$.  Now, we consider $\varphi':= (1+\alpha) \min (\varphi,\varphi_U)$ for $\alpha>0$. We have that $\varphi \leq -1/(1+\alpha)$ on some open neighborhood $W$ of $E\backslash U$ in the induced topology of $X \backslash U$ (that is $W=W' \backslash  U$ where $W'$ is an open set of $X$). So $\varphi'$ is less than $-1$ on the open set $V= U\cup W$, so it is less than $-1$ on some neighborhood of $E$. We let $\alpha$, $\varepsilon_1$, and $\varepsilon_2$ go to zero, and we deduce that $C(E)=C'(E)$. That concludes the first assertion.\\

\noindent For the second assertion, we show that $C(\cup E_i) \leq \lim C(E_i)$ (the other inequality is a consequence of the first item of proposition \ref{firstproperties}). For each $i$, let $\varphi_i \leq 0$,  with $\varphi_i = -1$ on $E\backslash H_{\varphi_i}$ where  $H_{\varphi_i}$ is pluripolar, such that $\|\varphi_i\|^2_* \leq C(E_i)+1/i$. Since $C\leq 1$, the sequence $(\varphi_i)$ is bounded in $W^*$ and we can extract a subsequence weakly converging in $W^*$. We apply lemma \ref{pointwise} for $a=-1$ and $A= E_i$, we obtain that the limit $ \varphi $ is equal to $-1$ on each $E_i$ (apart from some pluripolar set) and $\| \varphi \|_* \leq \lim \| \varphi_i\|_*$. That gives the proposition. $\Box$  \\

We also have the following description of $C$:
\begin{corollaire}
For all Borel set $E$, there exists an element $u_E\leq 0$ equal to $-1$ on $E\backslash H$ where $H$ is pluripolar with $\|u_E\|_*^2=C(E)$.
\end{corollaire}
\emph{Proof.} Take $(u_n)$ a sequence in $L(E)$ with $\|u_n\|_*^2\to C(E)$ and apply lemma \ref{pointwise} to $E$. $\Box$

\begin{Remark} \rm{It is not clear if such an extremal function is unique. It would also be interesting to know if the extremal function is qpsh or semi-continuous.}
\end{Remark}

\begin{Remark}\label{generalisation} \rm In the local case, we can in the same way define the capacity $C$. We show in the same way that it is a Choquet capacity for which the sets of zero capacity are exactly the pluripolar sets. The only difference is that we do not have that $C$ and ${\rm cap}_{BT} $ are comparable. But by \cite{AT1}, we can prove that they are locally comparable (it is the same argument as the one in the proof of proposition \ref{quasi-continuity3}). In particular, the functions in $W^*$ are locally quasi-continuous and the argument at the beginning of the proof of theorem \ref{quasi-continuity2} shows that if a function is locally quasi-continuous for a subadditive capacity then it is quasi-continuous. In  particular, the elements of $W^*$ are quasi-continous for the functional capacity $C$ which is a Choquet capacity.
\end{Remark}
Now, we consider  the set $M_\infty$ of positive Radon measures bounded for the norm $\|.\|_*$ on the space of smooth functions. That is $\mu\in M_\infty$ if $\mu$ is a positive Radon measure such that there exists a constant $A$ such that $|\int f d\mu|\leq A\|f\|_*$ for all $f$ smooth. We put on $M_\infty$ the operator norm $\|.\|'$ (that is the infimum of the $A$ above). 

Let $E$ be a Borel set, we define the capacity:
\begin{eqnarray*}
{\rm cap}'(E)&:=& \sup \Big\{\mu(E)^2, \ \mu \in M_\infty, \ \|\mu\|'\leq 1 \Big\}. 
\end{eqnarray*}
Observe that this set function defines an inner capacity. Recall the notation $L(E)=\{ \varphi \in W^*, \  \varphi \leq -1 \ \text{ a.e on some neighborhood of} \ E, \ \varphi \leq 0 \ \text{on} \ X\}$. We have the following proposition:
\begin{proposition}
Let $E$ be a Borel subset of $X$, then:
$${\rm cap}'(E)\leq C(E).$$
\end{proposition}
\emph{Proof.} It is sufficient to prove the inequality for $E$ compact. For any $\mu \in M_\infty$  with $\|\mu\|'\leq 1$ and $f$ in $L(E)$, we define $\langle\mu,f \rangle:=\int f d\mu$. Since $\mu$ is a positive Radon measure and $E$ is compact, we know that:
$$\mu(E)= \inf \Big\{ \langle \mu,-f \rangle, \ f\in L(E)\Big\}.$$  
By definition, for $\mu \in M_\infty$ and $f\in L(E)$, we have $0\leq \langle \mu,-f \rangle\leq \|f\|_* $. So we deduce from the previous inequality that $\mu(E)^2\leq C(E)$. Taking the supremum over all the measures $\mu$ in $M_\infty$ with $\|\mu\|'\leq1$ gives the result.   $\Box$ 

\begin{Remark} \rm{For Sobolev spaces, the functional capacity and its dual capacity coincide \cite{AH1}. It would be interesting to know if it is true here (this question is likely linked to the study of $u_E$).}
\end{Remark}
\begin{proposition}
Let $\mu \in M_\infty$, then $\mu$ does not charge pluripolar sets. For a Borel set $E$, ${\rm cap}'(E)=0$ if and only if $E$ is pluripolar.
\end{proposition}
\emph{Proof.} The first part is already in \cite{DS4} the second one follows from above. $\Box$

\begin{appendix} 
\section{Dirichlet space for a positive closed current}\label{appendix}
The results in this appendix are adaptated from \cite{Deny} and \cite{Oka1}.
\subsection{General setting}
Let $U$ be a bounded open subset of $\mathbb{C}^k$ with smooth boundary. Let $S$ be a positive closed current of bidimension $(1,1)$ (we do not assume that $S$ is of finite mass). We consider the quotient of the space of smooth forms by the kernel of the nonnegative bilinear form $ \langle u,v \rangle_S:=\text{Re} \int_U i\partial u \wedge \bar \partial v \wedge S$. Let $H_S$ be the completion of that quotient for the corresponding norm $\| . \|_S$. It is a Hilbert space. Define by $\mu:=i\partial\bar{\partial} \|z\|^2 \wedge S$ the trace measure of $S$. 

By lemma \ref{PoincaréSobolev} and the Cauchy-Schwarz inequality, there is a continuous inclusion from $H_S$ into $L^2(\mu)$ hence from $H_S$ into $L_{loc}^1(\mu)$ (for $V\Subset U$, then $\varphi \mapsto \varphi$ from $H_S(V)\to L^1(V)$ is continuous). Let $E\subset U$ be a Borel set. Define:
$$\mathcal{L}(E):=\Big\{ \varphi \in H_S, \ \varphi \leq 0, \ \varphi \leq -1 \ \text{a.e on a neighborhood of} \ E \Big\}.$$
 We define the capacity ${\rm cap}_S(E)$ by:
$${\rm cap}_S(E):=\inf \Big\{ \| \varphi\|_S ^2, \ \varphi \in \mathcal{L}(E) \Big\},$$
if $\mathcal{L}(E)$ is non empty, else we defined ${\rm cap}_S(E):=+\infty$.
\begin{proposition}
The following assertions hold:
\begin{enumerate}
\item For $E \subset F$, ${\rm cap}_S(E) \leq {\rm cap}_S(F)$.
\item For $(E_i)$ a sequence of Borel subsets of $U$,
$${\rm cap}_S(\cup E_i) \leq \sum {\rm cap}_S(E_i).$$
\item For $(K_i)$ a decreasing sequence of compacts, 
$$\lim {\rm cap}_S(K_i) = {\rm cap}_S( \cap K_i). $$
\item For $E_1 \subset E_{2} \subset \cdots $ an increasing sequence of Borel subsets of $U$, 
$${\rm cap}_S (\cup E_i) =\lim_{i \to \infty} {\rm cap}_S(E_i).$$
\end{enumerate}
Thus ${\rm cap}_S$ defines a Choquet capacity. 
\end{proposition}
\emph{Proof.} The first and third assertions are clear. For the second assertion, we can restrict ourselves to the case where $\sum {\rm cap}_S(E_i)$ is bounded. For all $i$, let $\varphi_i$ be in  $\mathcal{L}(E)$ such that ${\rm cap}_S(E_i) \geq \|\varphi_i\|_S^2 -2^{-i}\varepsilon$. Define $\psi_n:=\inf_{i\leq n} \varphi_i$, then we have that 
$$\|\psi_n\|_S^2 \leq \sum_{i\leq n} \|\varphi_i\|_S^2 \leq \sum_i {\rm cap}_S(E_i)+\varepsilon.$$
It is the same argument as in the proof of formula (\ref{extremum}) in section \ref{examples}: we prove the inequality $\|\inf_{i\leq n} \varphi_i\|_S^2 \leq \sum_{i\leq n} \|\varphi_i\|_S^2$ for smooth functions using a regularization of the functions $\max$ and $\min$ and we extend it to $H_S$ next. Taking a weakly converging subsequence in $H_S$, we have a non positive function $\psi \in H_S$  less than $-1$ in a neighborhood of $\cup E_i$ with  $\|\psi \|_S^2 \leq \sum_i {\rm cap}_S(E_i)+\varepsilon$.\\

There are several proofs for the fourth assertion. We follow the one in \cite{FZ1}. It relies on a geometric property of the norm $\| .\|_S$ that does not exist for $W^*$. Namely, following the proof of formula (\ref{extremum}) in section \ref{examples}, we have for $u$ and $v$ in $H_S$ that $\sup(u,v)$ and $\inf(u,v)$ are in $H_S$ with:
$$   \|\sup(u,v)\|_S^2+\| \inf(u,v)\|_S^2=\|u\|_S^2+\|v\|_S^2.$$
We restrict ourselves to the case where the sequence $({\rm cap}_S(E_i))_i$ is convergent. For all $i\geq 0$, let $u_i$ be such that $\|u_i\|_S^2\leq {\rm cap}_S(E_i)+\varepsilon_i$ with $u_i \in \mathcal{L}(E_i)$ and $\sum \varepsilon_i=\varepsilon$. Define $v_n=\inf_{i\leq n} v_i$. Observe that: $v_n=\inf(v_{n-1},u_n)$, $\|v_n\|_S^2\geq {\rm cap}_S(E_n)$ and $\|  \sup(v_{n-1},u_n) \|_S^2 \geq {\rm cap}_S(E_{n-1})$. So:
\begin{eqnarray*}
\| v_n \|_S^2+{\rm cap}_S(E_{n-1}) & \leq & \int_U i\partial \inf(v_{n-1},u_n) \wedge \bar \partial \inf(v_{n-1},u_n)\wedge S\\
                       &      &+ \int_U i\partial \sup(v_{n-1},u_n) \wedge \bar \partial \sup(v_{n-1},u_n)\wedge S\\
                       & \leq & \int_U i\partial v_{n-1} \wedge \bar \partial v_{n-1}\wedge S+\int_U i\partial u_n \wedge \bar \partial u_n \wedge S\\
                       & \leq & \| v_{n-1} \|_S^2+{\rm cap}_S(E_n)+\varepsilon_n
\end{eqnarray*}
Adding all these expressions from 1 to $n$ (with $\| h_1\|_S^2 \leq {\rm cap}_S(E_1)+\varepsilon_1$ for $n=1$), we get:
$$\|v_n \|_S^2\leq {\rm cap}_S(E_n)+\sum \varepsilon_i \leq \lim {\rm cap}_S(E_i) + \varepsilon.$$
We conclude by taking a weak limit of $(v_n)$. $\Box$

\begin{proposition}\label{lusin}
Let $u\in H_S$. Then $u$ is quasi-continuous: it admits a representative $\widetilde{u}$ such that for every $\varepsilon>0$ there exists an open subset $\Omega_\varepsilon$ with ${\rm cap}_S(\Omega_\varepsilon)\leq \varepsilon$ and $\widetilde{u}$ restricted to $U\backslash \Omega_\varepsilon$ is continuous.
\end{proposition}
\emph{Proof.} By definition of $H_S$, there is a sequence of smooth functions $(u_n)$ in $H_S$ converging to $u$ in $H_S$. We can suppose that $\|u_n-u_{n+1}\|_S^2\leq 2^{-n}$. Let $\widetilde{u}:=\lim u_m = u_m+\sum_{l> m} (u_{l+1}-u_l) $. 
Then, this series converges uniformly  on the closed set $\mathcal{E}_m:=\cap_{l>m} \{ |u_{l+1}-u_l|\leq l^{-2} \}$. So $\widetilde{u}$ is continuous there. Consider $v_m:=-\sum_{l>m} l^2 |u_{l+1}-u_l|$, it is less than $-1$ on  $\Omega_m:=U\backslash \mathcal{E}_m$ and it is in $H_S$. This gives the result for $m$ large enough. $\Box$
\begin{Remark} \rm{The previous proof is a simple consequence of the fact that smooth functions are dense in $H_S$. Proving a similar result in the case of $W^*$ is a main difficulty of this paper.} 
\end{Remark}
Recall that for $u\in H_S$, the current $\theta \mapsto \langle i\partial \bar{\partial}u \wedge S ,\theta \rangle:=-\langle u , \theta \rangle_S $ is well defined.
\begin{defi}
A function $\varphi$ is $S$-subharmonic if $i\partial\bar{\partial} \varphi \wedge S$ is a positive Radon measure on $U$.  
\end{defi} 
For $E\subset U$ with ${\rm cap}_S(E)\neq \infty$, the infimum in the definition of ${\rm cap}_S(E)$ is reached for a unique element $u_E$ in $H_S$. That is ${\rm cap}_S(E)=\|u_E\|_S^2$. Furthermore, $u_E$ is equal to $-1$ on the interior of $E$, it is $S$-subharmonic on $U$ and $S$-harmonic on any open set in $U\backslash E$. The following lemma is useful:
\begin{lemme}\label{tool}
Let $\varphi\in H_S$ and let $u$ be an $S$-subharmonic function in $H_S$ with $\varphi\leq u$. Then $\|u\|_S \leq \|\varphi\|_S$.
\end{lemme}
\emph{Proof.} Let $\theta:=u-\varphi \geq 0$. So $\langle i\partial \bar{\partial}u \wedge S ,\theta \rangle \geq0$. That means, by definition: 
$$- \langle u-\varphi,u \rangle_S \geq 0. $$
That is to say $\langle u,u \rangle_S \leq \langle \varphi,u \rangle_S$. By Cauchy-Schwarz inequality, that gives $\|u\|_S^2\leq \|u\|_S\|\varphi\|_S$ and the result follows. $\Box$

\subsection{Potentials}
Let $B$ be the set of function in $L^{\infty}(\mu)$ with compact support in $U$. For $f\in B$,  we have that $u \mapsto \int_U f u d\mu$ is a continuous linear form on $H_S$. Thus there exists an element $U^f$ such that for all $u\in H_S$,  $\int_U f u d\mu = \langle U^f,u \rangle_S$. Let $P$ be the closure of the elements $U^f$ with $f$ non negative; $P$ is a closed convex cone in $H_S$ (the cone of \emph{potentials}). Let $V$ be an open set in $U$, define $P(V)$ to be the closure in $H_S$ of the elements $U^f$ where $f\in B$ is non-negative with support in $V$. 

\begin{proposition}
Let $u \in H_S$, then the following properties are equivalent:
\begin{enumerate}
\item $u \in P(V)$.
\item For all $v \in H_S$, $v\geq 0$ on $V$, then $\langle v,u \rangle_S \geq 0$, or equivalently $\|v+u\|_S \geq \|u\|_S$.
\end{enumerate}
 In particular, $u_V$ is in $P(V)$.
\end{proposition}
\emph{Proof.} Let $u=U^f$ for some positive $f\in B$ with support in $V$ and let $v\in H_S$, $v\geq 0$ on $V$. Then $\langle v,u \rangle_S=\langle v,U^f \rangle_S=\int vf d\mu \geq 0$. By density, we conclude for the first implication.
 
 On the other hand, let $u$ be such that $\langle v,u \rangle_S \geq 0$ for all $v \in H_S$, $v\geq0$ on $V$. Let $u'$ be its projection on $P(V)$, it is characterized by $\langle u',u'\rangle_S=\langle u,u'\rangle_S$ and $\langle u',h\rangle_S\geq \langle u,h\rangle_S$ for all $h \in P(V)$. Thus, $\langle u',v\rangle_S\geq \langle u,v\rangle_S$, for $v=U^f$ with $f\in B$ positive, with support in $V$. That means $u' \geq u$ $\mu$-a.e in $V$. Now, we have $\|u'-u\|^2=\langle u',u'-u\rangle_S-\langle u,u'-u\rangle_S=-\langle u,u'-u\rangle_S\leq 0$. So $u'=u$ and $u$ is in $P(V)$. $\Box$\\

\noindent  Let $P^1(V)$ be the closure of the elements $U^f$ with $f\in B$ positive, with support in $V$, satisfying $\int f d\mu=1$. It is a closed convex set in $H$, empty if and only if $\mu(V)=0$. We give an alternate description of the capacity.

\begin{proposition}
Let $V$ be an open set, then we have the equality:
\begin{itemize}
\item ${\rm cap}_S(V)=0$ if $P^1(V)=\varnothing$;
\item ${\rm cap}_S(V)=\frac{1}{\inf_{u\in P^1(V)} \|u\|_S^2}$ if $P^1(V)\neq \varnothing$.
\end{itemize}
Furthermore in the case where  ${\rm cap}_S(V)$ is finite, $u_V=0$ if $\mu(V)=0$, else $u_V = -v/\|v\|_S^2$ where $v$ is the element of minimal norm in $P^1(V)$. 
\end{proposition}
\emph{Proof.} We assume $\mu(V)>0$, so $P^1(V)$ is not empty. Let $v$ be its element of minimal norm: $v$ satisfies $\langle v,U^f-v\rangle_S\geq 0$ for all $f\in B$ positive,  with support in $V$ satisfying $\int f d\mu=1$. So $v(x)\geq \|v\|_S^2$ $\mu$-a.e on $V$. 

If ${\rm cap}_S(V)=+\infty$ then $\mathcal{L}(V)$ is empty (i.e. there is no element $\geq 1$ on $V$). That means $\|v\|_S=0$ hence $(\inf\big\{\|u\|_S^2 , \ u\in P^1(V)\big\})^{-1} =+\infty$.

Assume, ${\rm cap}(V)$ is finite, then $u_V$ exists and satisfies for $f$ above $\langle -u_V,U^f\rangle_S=\int -u_Vf d\mu = 1$. So this inequality stands for $v$: $\langle v,-u_V\rangle_S = 1$. It implies $v\neq 0$ and so $(\inf\big\{\|u\|_S^2 , \ u\in P^1(V)\big\})^{-1}$ is finite. Set $w:=-v/ \|v\|_S^2$. We have $w\leq -1$ on $V$. From above, we deduce by Cauchy-Schwarz inequality that $\|w\|_S^2\leq \|u_V\|_S^2$. So $w=u_V$ by unicity of the element of minimum norm of $\mathcal{L}(V)$. $\Box$ \\

\noindent We can now prove a result that allows us to really define pointwise  values for the functions in $H_S$ up to a set of zero capacity.
\begin{theorem}
Let $u\in H_S$ be quasi-continuous, such that $u\leq 0$ $\mu$-a.e. Then, $u\leq 0$, quasi-everywhere.
\end{theorem}
\emph{Proof}
It is sufficient to show that the sets $E_\alpha=\{x| \ u(x)>\alpha>0\}$ are of zero capacity. Assume it is false, and choose $\alpha$ so that ${\rm cap}_S(E_\alpha)>0$. Let $V$ be an open set of capacity $<{\rm cap}_S(E_\alpha)$ such that $u$ is continuous on $U\backslash V$. The set $\Omega=E_\alpha \cup V$ is open because by continuity $E_\alpha \cap (U\backslash V):= \{x \in U\backslash V|u(x)>\alpha>0\}$ is an open set for the induced topology. Its capacity is stricly greater than ${\rm cap}_S(V)$.
From the previous proposition, there is a function $g\in B$ positive, satisfying $\int g d\mu=1$ with support in $\Omega$, such that $\|U^g\|^2 \leq \frac{1}{{\rm cap}_S(\Omega)}+\varepsilon$,  where $\varepsilon>0$ can be taken arbitrarily small. Furthermore:
$$\int_V g d \mu \leq -\int g u_V d\mu \leq \|u_V\|_S \|U^g\|_S\leq \sqrt{{\rm cap}_S(V)}\sqrt{\frac{1}{{\rm cap}_S(\Omega)}+\varepsilon},$$
and this last quantity is less than $1$ for $\varepsilon$ small enough. Thus: 
 $$\int_{E_\alpha \backslash V} g d\mu =\int_\Omega g d\mu -\int_V g d\mu >0,$$ 
a contradiction. $\Box$ \\

\noindent As in the proof of proposition \ref{lusin}, we obtain the following pointwise convergence result:
\begin{lemme}\label{lusin2}
Let $(u_n)$ be a strongly converging sequence in $H_S$, then we can extract a subsequence that converges outside a set of $S$-capacity zero. 
\end{lemme} 

\end{appendix}
\bibliography{biblio} 
\noindent Gabriel Vigny, Mathématiques - Bât. 425, UMR 8628,\\ 
 Université Paris-Sud, 91405 Orsay, France.  \\
\noindent Email: gabriel.vigny@math.u-psud.fr

\end{document}